\documentclass[a4paper]{amsart}

\usepackage{amsfonts,amssymb,amsthm,amsmath,amsxtra,amscd,verbatim,eucal}
 \def\cal#1{\mathcal{#1}}
 \usepackage[all]{xy}
 
\def\AA{\mathbb{A}}
\def\PP{\mathbb{P}}
\def\RR{\mathbb{R}}

\def\NN{\mathbb{N}}
\def\ZZ{\mathbb{Z}}

\let \cedilla =\c

\def\hmld{{\widehat{\mathrm{mld}}}}
\renewcommand{\a}{{\frak{a}}}
\renewcommand{\o}{{\mathcal O}}
\newcommand{\mld}{\mathrm{mld}}
\newcommand{\Hom}{\mathrm{Hom}}
\renewcommand{\j}{{\mathcal{J}}}
\newcommand{\cont}{{\mathrm{Cont}}}
\newcommand{\codim}{{\mathrm{codim}}}
\newcommand{\ord}{{\mathrm{ord}}}
\newcommand{\mult}{{\mathrm{mult}}}
\newcommand{\spec}{{\mathrm{Spec}}}

\newtheorem{Theorem}{Theorem}[section]
\newtheorem{Corollary}[Theorem]{Corollary}
\newtheorem{Corollary-Definition}[Theorem]{Corollary-Definition}
\newtheorem{Proposition}[Theorem]{Proposition}
\newtheorem{Lemma}[Theorem]{Lemma}

\newtheorem{conj}[Theorem]{Conjecture}

\theoremstyle{definition}
\newtheorem{defn}[Theorem]{Definition}
\newtheorem{exmp}[Theorem]{Example}
\newtheorem{rem}[Theorem]{Remark}
\newtheorem{say}[Theorem]{}

\newtheorem{Proposition-Definition}[Theorem]{Proposition-Definition}

\begin{document}

\title[]{Singularities with the highest Mather minimal log discrepancy}

\author{Shihoko Ishii and Ana J. Reguera}

\maketitle


\begin{abstract}
This paper characterizes  singularities with Mather minimal log discrepancies in the
highest unit interval, {\it i.e.,} the interval between $d-1$ and $d$, where $d$ is the dimension of the scheme.
The class of these singularities coincides with one of  the classes of (1) compound Du Val singularities, (2)  normal crossing double singularities, (3) pinch points, and
(4) pairs of non-singular varieties and  boundaries with  multiplicities less than or equal to $1$ at the point.
As a corollary, we also obtain one implication of an equivalence conjectured by Shokurov for the usual minimal log discrepancies. 
\end{abstract}

\section{Introduction}
\noindent
Let  $(X, B)$ be a pair 
consisting of a normal variety $X$ over an algebraically closed field $k$ of characteristic zero and an effective $\Bbb R$-divisor $B$ on $X$ such that $K_X+B$ is an $\Bbb R$-Cartier divisor.
The minimal log discrepancy $\mld (x, X, B)$ at a closed point $x\in X$ is defined for a pair  and plays an important role in birational geometry. 
On the other hand, we can also define Mather minimal log discrepancy $\hmld(x;X,{{\mathcal J}_X}B)$ with respect to the Jacobian ideal $\mathcal J_X$ of 
$X$ by using Mather discrepancy and the Jacobian ideal instead of the usual discrepancy (see [Is] and [DD]). 
Here  note that we need not  assume the $\Bbb R$-Cartier condition on $K_X+B$, and 
 $X$ can even be non-normal.
Mather minimal log discrepancy coincides with the usual discrepancy if $(X, x)$ is normal and locally a complete intersection.
We expect Mather minimal log discrepancy also to play  an important role in algebraic geometry, because  it sometimes has better properties than the usual minimal log discrepancy ([Is], [DD]).

Regarding the usual minimal log discrepancy,
Shokurov has proposed  the following conjectures:

\begin{conj} [{[Sh]}, Conjecture 2]
\label{sho1}
We have the inequality 
$$\mld (x;X, B)\leq \dim X,$$
where equality holds if and only if $(X,x)$ is non-singular and $B=0$ around $x$.
\end{conj}

Conjecture \ref{sho1} was proved for a non-degenerate hypersurface  case ([Am2]) and 
a three-dimensional Gorenstein case ([Ka], [Mar]); however,  it is still not proved in general.
But if one replaces $\mld$ by Mather minimal log discrepancy with respect to the Jacobian ideal, then the conjecture was proved  essentially in [Is], Corollary 3.15 (independently  proved also in [DD], Corollary 4.15).
Here we note that on a  variety $X$ (not necessarily normal), an effective $\RR$-Cartier divisor $B$ is 
defined as $B=\sum _{i=1}^sr_iB_i$ $(r_i\in \RR_{\geq 0})$, where $B_i$ is a subscheme on $X$  
defined by a principal ideal for $i=1,\ldots, s$. 

\begin{Proposition}[{[Is], Corollary 3.15; [DD], Corollary 4.15}]
\label{sh1-answer}
For an arbitrary variety $X$ and an effective $\RR$-Cartier divisor $B$ on $X$, we have the inequality 
$$\hmld(x;X, {{\mathcal J}_X}B)\leq \dim X,$$
where the equality holds if and only if $(X,x)$ is non-singular and $B=0$ around $x$.

\end{Proposition}

Shokurov has also proposed a conjecture (not published) as follows:

\begin{conj}
\label{sho2}
The inequality 
$$\dim X-1<\mld (x;X, B) $$ holds if and only if $(X,x)$ is non-singular and 
$\mult_x B<1$.
In this case, the minimal log discrepancy is computed by the exceptional divisor of the
first blow-up at $x$.
\end{conj}

The implication of the ``if" part of Conjecture \ref{sho2} for the two-dimensional case was proved 
by Vyacheslav Shokurov in an unpublished paper, and that for the three-dimensional case was proved by 
Florin Ambro [Am1]; however, this conjecture is not yet proved in general.
The main result of this paper is the following, which proves the Mather version of Conjecture  \ref{sho2}.

\begin{Theorem} A pair $(X, B)$ consisting of an arbitrary variety $X$ and an effective $\RR$-Cartier divisor $B$ on $X$ satisfies $$\dim X-1\leq {\hmld}(x;X,{{\mathcal J}_X}B) $$ 
if and only if
either 
\begin{enumerate}
\item[(i)]$B=0$, and $(X,x)$  is a normal crossing double singularity or a pinch point,
\item[(ii)]$B=0$, $\dim X\geq 2$ and $(X,x) $ is a compound Du Val singularity,   or
\item[(iii)]$(X,x)$ is non-singular and $0\leq \mult_xB\leq 1$.
\end{enumerate}

\noindent
In  cases $(\mathrm i)$ and $(\mathrm {ii})$, we have ${\hmld}(x;X,{{\mathcal J}_X})= \dim X -1$ and  
 in case $(\mathrm {iii})$, we have 
${\hmld}(x;X,{{\mathcal J}_X}B)=\mld(x; X, B)= \dim X -\mult_xB$, and the minimal  log discrepancy is computed by the exceptional divisor of the first blow-up at $x$.
\end{Theorem}

As a corollary, we obtain  the  ``if" part of Conjecture \ref{sho2} for the usual mld:
\begin{Corollary} The inequality $$\dim X-1<\mld (x;X, B)$$ holds
 if  $(X,x)$ is non-singular and 
$mult_x B<1$.
In this case the minimal log discrepancy is computed by the exceptional divisor of the
first blow-up at $x$.
\end{Corollary}
As a further corollary, we have the following  for usual mld.
\begin{Corollary}
\label{c.i.}
 Let $X$ be locally a complete intersection at a closed point $x\in X$.
Then, $\mld (x, X,\o_X)\leq \dim X$. Moreover, $\mld (x,X, \o_X)= \dim X$ if and only if $(X,x)$ is non-singular, and $\mld (x, X,\o_X)= \dim X-1$ if and only if $(X, x)$ is either a cDV singularity, a normal crossing double singularity, or a pinch point.
\end{Corollary}



\noindent
{\bf Acknowledgement:} The authors are grateful to Vyacheslav Shokurov, Florin Ambro and Yuri Prokhorov for 
 useful discussions and information about the usual minimal log discrepancies.

\section{Preliminaries on arcs and Mather discrepancy}

\begin{say} Let $k$ be a field and let $X$ be a $k$-scheme. Given $m \in
\ZZ_{\geq 0}$ and a field extension $k \subseteq K$, a $K$-arc of
order $m$ (resp. a $K$-arc) on $X$ is a $k$-morphism $\text{Spec }
K[[t]]/ (t)^{m+1} \rightarrow X$ (resp. $\text{Spec } K[[t]]
\rightarrow X$). There exists a $k$-scheme $X_m$ of finite type over $k$,
called the {\it space of arcs of order} $m$ of $X$, whose
$K$-rational points are the $K$-arcs of order $m$ on $X$, for any
$K \supseteq k$. There are natural affine morphisms $j_{m',m}:
X_{m'} \rightarrow X_m$ for $m' > m$. The projective limit
$X_\infty:= \lim_{\leftarrow} X_m$ is a 
$k$-scheme (not of finite type), called the {\it space of arcs of} $X$, whose $K$-rational
points are the $K$-arcs on $X$. We denote the natural morphisms by $j_m: X_\infty
\rightarrow X_m$, $m \geq 0$. \\

The scheme $X_\infty$ satisfies the following representability
property: for every $k$-algebra $A$, we have a natural isomorphism
$$
\Hom_k(\text{Spec } A, X_\infty) \cong \Hom_k(\text{Spec } A[[t]],
X) .
$$
Given $P \in X_\infty$  with residue field $\kappa(P)$, let us
denote by $h_P$ the $\kappa(P)$-arc on $X$  corresponding to the
$\kappa(P)$-rational point of $X_\infty$ by the previous
isomorphism. The image in $X$ of the closed point of $\text{Spec}
\ \kappa(P)[[t]]$, or equivalently, the image of $P$ by
$j_0:X_\infty \rightarrow X=X_0$ is called the {\it center} of
$P$. Then, $h_P$ induces a morphism of $k$-algebras $h_P^* :
{\cal O}_{X, j_0(P)} \rightarrow \kappa(P) [[t]]$.
	

For any proper closed subset $W$ of $X$, let $X_\infty^W:=j_0^{-1}(W)$
(resp. $X_m^W:=j_{m,0}^{-1}(W)$ for $m \geq 0$) be the closed
subset of $X_\infty$ (resp. of $X_m$) consisting of the arcs
(resp. arcs of order $m$) on $X$, whose center lies on $W$. If
$x_0$ is a fixed point in $X$, we simplify the previous notation
setting $X_m^0:=X_m^{x_0}$ for $m \geq 0$, and
$X_\infty^0:=X_\infty^{x_0}$ when there is no risk of confusion.
In particular, if $R$ is a local ring with  maximal ideal $M$, then
$(\text{Spec } R)_m^0:=(\text{Spec } R)_m^M$, for $m \geq 0$, and
$(\text{Spec } R)_\infty^0:=(\text{Spec } R)_\infty^M$.
\end{say}

\begin{say} 
\label{notation}
The space of arcs of the affine space $\AA^N_k= \text{Spec }
k[x_1, \ldots, x_N]$ is
$$
(\AA^N_k)_\infty= \text{Spec } k[{\underline {X}}_0, {\underline
{X}}_1, \ldots , {\underline {X}}_n, \ldots ],
$$
where for $n \geq 0$, ${\underline {X}}_n= ({X}_{1,n}, \ldots ,
{X}_{N,n})$ is an $N$-uple of variables. The space of arcs of
$\text{Spec } k[[x_1, \ldots, x_N]]$ is
$$ (\text{Spec } k[[x_1, \ldots, x_N]])_\infty = \text{Spec }
k[[{\underline {X}}_0]] [{\underline {X}}_1, \ldots , {\underline
{X}}_n, \ldots ] .
$$
For any $f \in k[x_1, \ldots , x_N]$ (resp. $f \in k[[x_1, \ldots
, x_N]]$) let $ \sum_{n=0}^\infty F_n \ t^n$ be the Taylor
expansion of $f(\sum_n {\underline {X}}_n \ t^n)$, hence $F_n \in
k[{\underline {X}}_0, \ldots, {\underline {X}}_n]$ (resp. $F_n \in
k[[{\underline {X}}_0]] [{\underline {X}}_1, \ldots, {\underline
{X}}_n]$). If $X \subseteq \AA^N_k$ is an affine variety, and $I_X
\subset k[x_1, \ldots , x_N]$ is the ideal defining $X$ in
$\AA^N_k$, then we have
$$
\aligned {X_m}= \text{Spec } k[{\underline {X}}_0, {\underline
{X}}_1, \ldots , {\underline {X}}_m] \ / \ (\{{F}_n\}_{0 \leq n
\leq m, f \in I_X}) \ \ \ \ \text{for } m \geq 0 \\
{X_\infty}= \text{Spec } k[{\underline {X}}_0, {\underline {X}}_1,
\ldots , {\underline {X}}_n, \ldots ] \ / \ (\{{F}_n\}_{n \geq 0,
f \in I_X}). \ \ \ \ \ \ \ \ \ \ \ \ \ \ \ \
\endaligned
$$
Suppose that ${\underline 0} \in X$. For $n \geq 0$, let $F_n^0$
denote the image of $F_n$ by the projection $k[{\underline {X}}_0,
{\underline {X}}_1, \ldots , {\underline {X}}_n] \rightarrow
k[{\underline {X}}_1, \ldots , {\underline {X}}_n]$, sending
$X_{i,0}$ to $0$ for $1 \leq i \leq N$. Note that $F_0^0=0$ for
$f \in I_X$.
Then
\begin{equation}
\label{1} X_m^0={X_m^{\underline 0}}= \text{Spec } k[{\underline
{X}}_1, \ldots , {\underline {X}}_m] \ / \ (\{{F}^0_n\}_{1 \leq n
\leq m,
f \in I_X}) \ \ \ \ \text{for } m \geq 0 
\end{equation}
$$\ \ \ \ 
X_\infty^0={X_\infty^{\underline 0}}= \text{Spec } k[{\underline
{X}}_1, \ldots , {\underline {X}}_m, \ldots ] \ / \
(\{{F}^0_n\}_{n \geq 0, f \in I_X}) . \ \ \ \ \ \ \ \ \ \ \ \ \ \
\
$$
Analogously, if $X \subset \text{Spec } k[[x_1, \ldots , x_N]]$
then
$$
\aligned {X_m}= \text{Spec } k[[{\underline {X}}_0]] [{\underline
{X}}_1, \ldots , {\underline {X}}_m] \ / \ (\{{F}_n\}_{0 \leq n
\leq m, f \in I_X}) \ \ \ \ \text{for } m \geq 0 \\
{X_\infty}= \text{Spec } k[[{\underline {X}}_0]] [{\underline
{X}}_1, \ldots , {\underline {X}}_n, \ldots ] \ / \ (\{{F}_n\}_{n
\geq 0,
f \in I_X}) \ \ \ \  \ \ \ \ \ \ \ \ \ \ \ \ \\
\endaligned
$$
and
\begin{equation}
\label{2}
 {X_m^{0}}= \text{Spec } k[{\underline {X}}_1, \ldots ,
{\underline {X}}_m] \ / \ (\{{F}^0_n\}_{1 \leq n \leq m, f \in
I_X}) \ \ \ \ \text{for } m \geq 0 \ \ \ \ \ \ \ \\
\end{equation}
$$\ \ \ \ \
{X_\infty^{0}}= \text{Spec } k[{\underline {X}}_1, \ldots ,
{\underline {X}}_m, \ldots ] \ / \ (\{{F}^0_n\}_{n \geq 0, f \in
I_X}) . \ \ \ \ \ \ \ \ \ \ \ \ \ \ \ \ \ \ \ \ \ \
$$
Here $F_n^0$ is the image of $F_n$ by the projection
$k[[{\underline {X}}_0]] [{\underline {X}}_1, \ldots , {\underline
{X}}_n] \rightarrow k[{\underline {X}}_1, \ldots ,
{\underline {X}}_n]$, sending $X_{i,0}$ to $0$, for $1 \leq i \leq N$.  
\end{say}

\begin{say}
\label{2.3}
 Given a germ of an algebraic variety $(X,x_0)$, {\it i.e.,} $X$ is
 a reduced separated $k$-scheme of finite type and $x_0$ is
 a closed point of $X$, recall that $X_m^0:=X_m^{x_0}$ for $m \geq
 0$ and $X_\infty^0:=X_\infty^{x_0}$. Let $R={\cal O}_{X,x_0}$
 and  $M$ be its maximal ideal. Note that
 we have
 \begin{equation}
 \label{3}
  X_m^0 \cong (\text{Spec } R)_m^0  \ \ \text{for }
m \geq 0 \ \ \text{ and } \ \ X_\infty^0 \cong (\text{Spec }
R)_\infty^0 .
\end{equation}

Let $\widehat{R}$  denote the $M$-adic completion of $R$ (the
notation $\widehat{{\cal O}_{X,x_0}}$ will also be used in text).
Then, every $K$-arc of order $m$ on $\text{Spec } R$ centered at
$M$ (resp. every $K$-arc on $\text{Spec } R$ centered at $M$)
extends in a unique way to a $k$-morphism \linebreak $\text{Spec }
K[[t]] \ / \ (t)^{m+1} \rightarrow \text{Spec } \widehat{R}$
(resp. $\text{Spec } K[[t]] \rightarrow \text{Spec }
\widehat{R}$), that is, a $K$-arc of order $m$ (resp. $K$-arc) on
$\text{Spec } \widehat{R}$, and it follows that
\begin{equation}
\label{4}
(\text{Spec } R)_m^0
\cong (\text{Spec } \widehat{R})_m^0 \ \ \text{for } m \geq 0 \ \
\ \text{ and } \ \ \
(\text{Spec } R)_\infty^0
\cong
(\text{Spec } \widehat{R})_\infty^0 .
\end{equation}

In fact, we may suppose that $X \subseteq \AA^N_k$ is affine, and
then, applying (\ref{3}),  equalities (\ref{1}) for $X \subseteq \AA^N_k$
and  equalities (\ref{2}) for $\text{Spec } \widehat{R} \subseteq
\text{Spec } k[[x_1, \ldots , x_N ]]$, an explicit description of
the
isomorphisms in (\ref{4}) follows. 
\end{say}

\begin{say}
\label{codim}
 Henceforth, $k$ will be an algebraically closed field of
characteristic zero, and $X$ an algebraic variety over $k$ of
dimension $d$.  Given a resolution of  the singularities $\pi: Y
\rightarrow X$ and a
 prime divisor $E$ on $Y$ contained in the exceptional locus of $\pi$, recall that $Y_\infty^{E}$ is the inverse
image of $E$ by the projection $Y_\infty \rightarrow Y$ and let
$N_E$ be the closure of its image $\pi_\infty (Y_\infty^{E})$ by
$\pi_\infty: Y_\infty \rightarrow X_\infty$, which is an
irreducible subset of $X_\infty^\text{Sing X}$. The generic point
$P_E$ of $N_E$ is a stable point of $X_\infty$ (see [Re], 3.1).
Therefore, $\dim {\cal O}_{\overline{j_m(X_\infty)}, j_m(P_E)}$ is
constant for $m\gg 0$ and we have
\begin{equation}
\label{5}
\dim {\cal O}_{X_\infty, P_E} \leq \sup_m \dim {\cal
O}_{\overline{j_m(X_\infty)}, j_m(P_E)} < \infty 
\end{equation}

Although the first inequality may be strict (see Example \ref{ana}
below), the constant $\dim {\cal O}_{\overline{j_m(X_\infty)},
j_m(P_E)}$, for $m\gg0$, is called the codimension of $N_E$ in
$X_\infty$ (see [EM], sec. 5). In [DEI], the value of this constant
has been described in terms of the relative Mather canonical
divisor $\widehat{K}_{Y/X}$. A review of this
concept is given below.
\end{say}

\begin{say}
\label{nash}
 If $\pi: Y \rightarrow X$ is a resolution of the singularities
dominating the Nash blowing-up of $X$ (for the definition of Nash blowing-up, see, for example 
[DEI], Definition 1.1), then the image of the
canonical homomorphism $ d \pi: \pi^*(\wedge^d \Omega_X)
\rightarrow \wedge^d \Omega_Y$ is an invertible sheaf. More
precisely, there exists an effective divisor $\widehat{K}_{Y/X}$
with support in the exceptional locus of $\pi$ such that
$$
d \pi (\pi^* (\wedge^d \Omega_X))= {\cal O}_Y(- \widehat{K}_{Y/X})
\ \wedge^d \Omega_Y.
$$
The divisor $\widehat{K}_{Y/X}$ is called the {\it Mather
discrepancy divisor}.\\

For any prime divisor $E$ on $Y$, let
$$
{\widehat k}_E:=\text{ord}_E(\widehat{K}_{Y/X}) .
$$
Here ${\widehat k}_E$ is an integer because $\widehat{K}_{Y/X}$ is a divisor on $Y$.
Note that ${\widehat k}_E \neq 0$ implies that $E$ is contained in
the exceptional locus of $\pi$ and that ${\widehat k}_E$ 
depends only on the divisorial valuation $\nu_E$ defined by $E$.
That
is, if $\pi: Y' \rightarrow X$ is another resolution of
singularities dominating the Nash blowing-up of $X$ such that the
center of $\nu_E$ on $Y'$ is a
divisor $E'$, then 
$\text{ord}_E(\widehat{K}_{Y/X})
=\text{ord}_{E'}(\widehat{K}_{Y'/X})$. 
In general, for every prime divisor $E$ over $X$, {\it i.e.,} a prime divisor on a normal variety $Y$, which is 
proper birational over $X$, we can define ${\widehat k}_E$ because there exists a resolution
 $Y'$ of the singularities that dominates both the Nash blowing-up and $Y$, and $E$ appears  
on $Y'$ also.

Then we have
\begin{equation}
\label{6}
\sup_m \dim {\cal O}_{\overline{j_m(X_\infty)}, j_m(P_E)}=
\widehat{k}_E+1
\end{equation}
([DEI], Theorem 3.9).\\
\end{say}

For each divisorial valuation $\nu=\nu_E$, where $E$ is a prime
divisor over $X$, let $c_X(\nu_E)$ denote the center of $\nu_E$
on $X$. Given non-zero ideal sheaves  in ${\cal O}_X$ with  non-negative powers
and a closed subset $W$ of $X$,  {\it Mather minimal
log discrepancy} along $W$ is defined as
follows:

\begin{defn} Let $X$ be a variety over $k$.
 Given  non-zero ideal sheaves $\a_1,\ldots , \a_l$ of ${\cal O}_X$, non-negative real numbers 
 $s_1,\ldots, s_l$, and a proper 
closed subset $W$ of $X$,  {\it Mather minimal log discrepancy
of } $(X, \a_1^{s_1},\ldots, \a_l^{s_l})$ along $W$ is defined as follows:

\vskip.5truecm
\begin{equation}
\label{7}
\widehat{\mld}(W;X, \a_1^{s_1}\cdots\a_l^{s_l})\ \ \ \ \ \ \ \ \ \ \ \ \ \ \ \ \ \ \ \ \ \ \ \ \ \ \ \ \ \ \ \ \ \ \ \ \ \ \ \
\end{equation}
$$:= {\inf} \ \{ \widehat{k}_E -
\sum_{i=1}^{l}s_i\mbox{ord}_E(\a_i) +1 \mid \ \nu_E \mbox{ divisorial
valuation}, c_X(\nu_E) \subseteq W\},
$$
if $\dim X \geq 2$ or $\dim X=1$ and the infimum on the right hand
side is non-negative; otherwise, $\widehat{\mld}(W;X,
\a_1^{s_1}\cdots\a_l^{s_l}):=-\infty$. 
Here  $\widehat{k}_E=\ord _E\widehat{K}_{Y/X} $ for a resolution $Y\to X$ on which $E$ appears.
A remark on this definition: if $\dim X
\geq 2$, then the infimum in (\ref{7}) is negative if and only if it is
equal to $-\infty$ ([Is], Remark 3.4).\\

Let $X$ be a  variety and $B=\sum_{j=1}^r b_jB_j$  an effective $\Bbb R$-Cartier divisor, 
{\it i.e.,} $b_j\in \Bbb{R}_{\geq 0}$ and $B_j$ is a subscheme defined by a principal ideal on $X$ $ (j=1,\ldots, r)$.
For non-zero ideal sheaves $\a_1,\ldots , \a_l$ of ${\cal O}_X$, positive real numbers 
 $s_1,\ldots, s_l$, and an effective $\Bbb R$-divisor $B$, we define Mather minimal log discrepancy for the mixture  of ideals and a divisor $(\a_1^{s_1},\ldots, \a_l^{s_l}, B)$ as follows
 $$
\widehat{\mld}(W;X, \a_1^{s_1}\cdots\a_l^{s_l}\cdot B):=
 \widehat{\mld}(W;X, \a_1^{s_1}\cdots \a_l^{s_l}\cdot\o_X(-B_1)^{b_1}\cdots \o_X(-B_r)^{b_r}).$$
\end{defn}

\begin{say}
\label{usualmld}
 If $X$ is a normal affine Gorenstein variety and $\pi: Y
\rightarrow X$ is a resolution of the singularities dominating the
Nash blowing-up of $X$, then we have
\begin{equation}
\label{8}
{\cal O}_Y(- \widehat{K}_{Y/X})= \pi^*(I_{Z}) \otimes {\cal O}_Y(-
K_{Y/X}), 
\end{equation}
where $K_{Y/X}$ is the unique divisor with support in the
exceptional locus of $\pi$, which is linearly equivalent to $K_Y -
\pi^*(K_X)$, and $I_{Z}$ is the ideal defining the first Nash
subscheme of $X$ ([EM], Appendix).\\

In particular, if $X$ is normal and  a complete intersection, then
$I_{Z}= {\cal J}_X$ is the Jacobian ideal of $X$, and from (\ref{8}), it
follows that for any non-zero ideal sheaf $\mathfrak{a}$ of ${\cal O}_X$
and any proper closed subset $W$ of $X$, we have
\begin{equation}
\label{9}
\widehat{\mld}(W;X, \a_1^{s_1}\cdots\a_l^{s_l}{\cal J}_X) = {\mld}(W;X,
\a_1^{s_1}\cdots\a_l^{s_l}), 
\end{equation}
where recalling that $\mld(W;X,  \a_1^{s_1}\cdots\a_l^{s_l})$ is the minimal log
discrepancy of $(X, \a_1^{s_1}\cdots\a_l^{s_l})$ along $W$ defined as follows:

\vskip.5truecm
${\mld}(W;X, \a_1^{s_1}\cdots\a_l^{s_l})$
$$:= \text{inf} \ \{ k_E -
\sum_{i=1}^{l}s_i\mbox{ord}_E(\a_i) +1 \mid \ \nu_E \text{ divisorial
valuation}, c_X(\nu_E) \subseteq W\}
$$
if $\dim X \geq 2$ or $\dim X=1$ and the infimum on the right hand
side is non-negative; otherwise, $\widehat{\mld}(W;X,
\a_1^{s_1}\cdots\a_l^{s_l}):=-\infty$. Here $k_E:= \text{ord}_E K_{Y /X}$,
where $Y \rightarrow X$ is a desingularization and $E$ is the
center of $\nu_E$ on $Y$ and only depends on
the divisorial valuation $\nu_E$.
\end{say}

\begin{exmp}
\label{ana} Let $X$ be the hypersurface $x_1^2+x_2^2+x_3^2=0$ in
$\AA^3_k$, which has an ${\bf A}_1$-singularity at $\underline 0$.
The exceptional locus of its minimal desingularization $\pi:Y
\rightarrow X$ consists of a unique irreducible curve $E$. We have
$\dim {\cal O}_{X_\infty, P_E}=1$ ([Re] Corollaries 5.12 and
5.15). On the other hand, $X$ has a canonical singularity at
$\underline 0$, hence $k_E=0$ and, by (\ref{8}), $\widehat{k}_E=1$. From
(\ref{6}), we conclude that
$$
1=\dim {\cal O}_{X_\infty, P_E} < \sup_m \dim {\cal
O}_{\overline{j_m(X_\infty)}, j_m(P_E)}=\widehat{k}_E + 1 =2
$$
That is, in this example, the first inequality in (\ref{5}) is strict.
\end{exmp}



\begin{say} {\it Inversion of adjunction} ([Is], Proposition  3.10, [DD], Theorem 4.10): Let $X$ be
an algebraic variety, $A$ a non-singular variety containing
$X$  as a closed subvariety of codimension $c$, and $W$ a proper closed
subset of $X$. Let $\widetilde{\mathfrak{a}} \subseteq {\cal
O}_A$ be an ideal sheaf such that
$\mathfrak{a}:=\widetilde{\mathfrak{a}} {\cal O}_{X}$ is a non-zero
ideal sheaf of ${\cal O}_{X}$, and let $I_X \subseteq {\cal
O}_A$ be the ideal sheaf defining $X$. Then
\begin{equation}
\label{inv}
\widehat{\mld}(W;X, \mathfrak{a}{\cal J}_X)=\widehat{\mld}(W;{A}, \widetilde{\mathfrak{a}}I_X^c) .
\end{equation}

This result is a generalization of an analogous result for minimal
log discrepancies proved in [EM], Theorem 8.1. 

Consider   equality (\ref{inv}) for trivial $\a, \tilde\a$ and $W=\{x_0\} $ for a closed point $x_0$. By [Is], Proposition 3.7, the right hand side is represented as follows: 
 \begin{equation}
 \label{repr}
 \hmld(x_0;A,I_X^c)=\inf_m\{\codim (\cont^{\geq m}(I_X)\cap \cont^{\geq 1}(M_{x_0}), A_\infty)-cm\},
 \end{equation}
 where the codimension, as in the sense of [EM] (see \ref{codim}), is the minimal value of the  $\{\codim (C_i, A_\infty)\}$ of the irreducible components 
 $C_i$ in $A_\infty$. 
 Here  if $P_i$ is the generic point of the component $C_i$, then 
 $\codim (C_i, A_\infty)$ is defined as  $\sup_m \dim \o_{\overline{j_m(A_\infty)},j_m(P_i)}$ 
and it coincides with the constant
 $$\dim \o_{j_m(A_\infty),j_m(P_i)}=\dim\o_{A_m,j_m(P_i)} \ \ {\mbox {for}} \ m\gg 0.$$
 Here by using the canonical projection $j_{r}:A_\infty\to A_{r}$,  the ``contact loci" are 
 represented as  follows:
 $$\cont^{\geq m}(I_X):=
\{P\in A_\infty \mid \ord_t h_P^*(I_X)\geq m\}=j_{m-1}^{-1}(X_{m-1}),$$ 
$$\cont^{\geq 1}(M_{x_0}):=\{P\in A_\infty \mid \ord_t h_P^*(M_{x_0})\geq 1\}=j_{0}^{-1}(x_0).$$
Therefore, (\ref{repr}) turns out to be
\begin{equation}
\label{repr2}
 \hmld(x_0;A,I_X^c)=\inf_m\{ \codim(X_{m-1}\cap j_{m-1,0}^{-1}(x_0), A_{m-1})-cm\}.
 \end{equation}
 From (\ref{inv}), (\ref{repr2}), and also replacing $m-1$ with $m$, we obtain
 
\begin{equation}
\label{formula}
\widehat{\mld}(x_0;X, {\cal J}_X)  =\inf_m \ \{ \ (m+1) d  - \dim
X^{0}_m \ \},
\end{equation}
where $d=\dim X$.
\end{say}

\begin{Proposition} 
\label{basic}
For an arbitrary variety $X$ of dimension $d$ and an effective $\RR$-Cartier divisor $B$ on $X$, we have the inequality 
$$\hmld(x;X, {{\mathcal J}_X}B)\leq d,$$
where the equality holds if and only if $B=0$ around $x$  and $(X,x)$ is non-singular.
\end{Proposition}

{\it Proof:}
In [Is], Corollary 3.15, it was proved that 
$$\hmld(x;X, {{\mathcal J}_X})\leq d$$
always holds
and the equality holds if and only if  $(X,x)$ is non-singular.
Therefore, it is sufficient to prove that  
$\hmld(x;X, {{\mathcal J}_X}B)< d$ for any  $B$ non-zero around $x$.
And this inequality holds because $\hmld(x;X, {{\mathcal J}_X}B)< \hmld(x;X, {{\mathcal J}_X})$ holds for any  $B$ non-zero around $x$. $\Box$

\section{Characterization of top singularities}

\begin{defn} For $d\geq 1$, we say that a $d$-dimensional variety $X$ has a {\it top singularity} at $x_0$ or that
$(X, x_0)$ is a top singularity, if $\widehat{\mld}(x_0;X, {\cal J}_X)=d-1$.\\
\end{defn}

\begin{defn}
Let $(X, {x_0})$ be a germ of a  hypersurface in $\AA^{d+1}_k$ with $d\geq 2$.
A singularity $(X, {x_0})$ is called a {\it compound Du Val singularity} or cDV singularity
if (i) it is Du Val in the case of $d=2$, (ii) its general hyperplane section is a Du Val singularity in the case of
$d=3$, (iii)  its general hyperplane section is a cDV singularity in the case of $d> 3$.
\end{defn}

In the definition of cDV singularities, we assume the generality of hyperplane sections,
but it is not necessary if one assumes that the point is originally singular.
The following is well known, but we present a proof here for  readers' convenience.

\begin{Lemma}
\label{special}
 Assume that $(X, {x_0})$ is a germ of a $d$-dimensional hypersurface singularity of multiplicity $2$ with $d\geq 3$,
 then $(X, {x_0})$ is a compound Du Val singularity if and only if there exist $(d-2)$ hyperplanes $H_1, \ldots, H_{d-2}$ such that $(X\cap H_1\cap\cdots\cap H_{d-2}, x_0)$
 is a Du Val singularity.
\end{Lemma}

{\it Proof:}
The ``only if" part is trivial and we will show the ``if" part.
Because the statement is local, we may assume that $(X, x_0)\subset (W, \underline 0)$, where 
$(W,\underline 0)$ is  an open neighborhood of $x_0=\underline{0}$ in $ (\AA^{N}_k, \underline 0)$, $(N=d+1)$.
Let $H_i=H_{{\bf a}_i}:=\{(x_1, \ldots, x_N)\in W\mid \sum_{j=1}^{N}a_i^jx_j=0\}$,
where ${\bf a}_i=(a_i^1,\ldots, a_i^{N})\in \AA^{N}_k$.
Define the family $\mathcal S$
$$\begin{array}{ccc}
 \mathcal {S}& \subset &  X\times \AA^{(d-2)N}_k\\ 
 &&\\
 &\rho\searrow&\downarrow\tilde\rho\\
 && \AA_k^{(d-2)N}\\
\end{array}$$
so that for $\underline{\lambda}=({\underline\lambda}_1, \ldots, {\underline\lambda}_{d-2})\in \AA_k^{(d-2)N}$
$({\underline\lambda}_i=(\lambda_i^1,\ldots, \lambda_i^N)\in \AA_k^{N})$,
$${\mathcal S}_{\underline\lambda}:=\rho^{-1}({\underline\lambda})=X\cap\left\{ \sum_{j=1}^{N}\lambda_i^jx_j=0\ \mbox{for}\ i=1,\ldots, d-2\right\}.$$
The space $\mathcal S$ is a successive $(d-2)$ hyperplane cut of 
$X\times \AA_k^{(d-2)N}$  around $({\underline 0},{\bf a})$ $({\bf a}=
({\bf a}_1,\ldots, {\bf a}_{d-2}))$ and the $\codim _{\underline 0}({\mathcal S}_{\bf a}, X)=d-2$, because ${\mathcal S}_{\bf a}$ is a surface with a Du Val singularity at $\underline 0$.
Therefore, $\codim _{(\underline 0,{\bf a})}({\mathcal S}, X\times\AA^{(d-2)N})=d-2$. 
Then, $\rho$ is flat around $({\underline 0},{\bf a})\in {\mathcal S}$ by a general theory
(see for example, [M] p. 177, Corollary).
By replacing  $W$, therefore also $X$, by a sufficiently small neighborhood,
there is an open neighborhood $U\subset \AA_k^{(d-2)N}$ of $\bf a$ such that 
${\mathcal S}|_U\to U$ is flat.
Then, the family  ${\mathcal S}|_U\to U$ is a flat family of surface singularities.
As $(X, x_0)$ is singular, $({\mathcal S}_{\underline\lambda}, \underline 0)$ is also singular for every point $\underline\lambda \in U$.
As is well known, Du Val singularities deform only to Du Val singularities (for example, see 
[B] or [KS]), 
and there is a neighborhood $U_0\subset U$ of $\bf a$ such that 
the singularity 
$({\mathcal S}_{\underline \lambda}, \underline 0)$ is  a Du Val singularity for $\underline \lambda\in U_0$. This proves that the general $(d-2)$ hyperplane cut of $(X,x_0)$ is Du Val. 
$\Box$

\begin{Corollary}
\label{3.4}
Let $(X, x_0)$ be a germ of a $d$-dimensional variety and $\widehat R$ be the $M$-adic completion of 
$R=\o_{X, x_0}$. Then, $(X, x_0)$ is a compound Du Val singularity if and only if there exists 
$g_1,\ldots, g_{d-2}\in \widehat R$ such that $\spec \widehat R/(g_1,\ldots, g_{d-2}) $ is a Du Val 
singularity.
\end{Corollary}

{\it Proof:} The necessary condition is trivial. For the
sufficient condition, let $(X,x_0) \subset (\AA^{d+1}_k,
{\underline 0})$ and let $x_1, \ldots, x_{d+1}$ be   a system of
coordinates in $\AA^{d+1}_k$. Note  that the condition of
$(X,x_0)$ being a cDV singularity does not depend on the choice of
the system of coordinates in $\AA_k^{d+1}$ (this follows from a
similar argument to that in the proof of Lemma \ref{special}). Let $g_1,
\ldots, g_{d-2} \in \widehat{R}$ be such that $R':= \widehat{R} /
(g_1, \ldots, g_{d-2})$ has a Du Val singularity. Then $\{g_1,
\ldots, g_{d-2}\}$ is a regular sequence consisting of elements of
multiplicity $1$; hence, for fixed $n\gg 0$, after a change of
affine coordinates in $\AA^{d+1}_k$, we may suppose that $g_i=x_i
\ \text{mod } M^n$, for $1 \leq i \leq d-2$, where $M$ is the
maximal ideal of $\widehat{R}$. From this it follows that
$\widehat{R}/(x_1, \ldots, x_{d-2})$ and  $R/(x_1, \ldots,
x_{d-2})$ have a Du Val singularity. Thus, $(X,x_0)$
is a cDV singularity by Lemma \ref{special}. $\Box$

\begin{Lemma}
\label{ci}
Let $X$ be a $d$-dimensional variety and let $X' \subset X$ be a
$(d-c)$-dimensional subvariety that is defined as the zero locus
of $c$ elements of ${\cal O}_X$. Let $x_0$ be a closed point in
$X'$. If $(X', x_0)$ is a top singularity,  then $(X, x_0)$ is also a
top singularity.

Moreover, given a germ of a $d$-dimensional variety $(X,x_0)$, let
$\widehat{R}$ be the $M$-adic completion of $R={\cal O}_{X,x_0}$
(see notation in \ref{2.3}). Suppose that there exist $g_1, \ldots, g_c
\in \widehat{R}$ such that the ring $R'=\widehat{R}/(g_1, \ldots,
g_c )$ has dimension $d-c$ and
$$
\inf_m \ \{ \ (m+1) (d-c)  - \dim (\text{Spec } R')^{0}_m \ \} =
d-c-1.
$$
Then $(X,x_0)$ is a top singularity. \\
%
\end{Lemma}

{\it Proof:} 
By (\ref{4}) and (\ref{formula}), it is sufficient to prove the second assertion.
We may suppose that $\widehat{R}=k[[x_1,\ldots,x_N]]/I$ for an ideal $I\subset  k[[x_1,\ldots,x_N]]$.
Let $g_1, \ldots , g_c \in k[[x_1,
\ldots , x_N]]$ be as in the lemma, and set $X'=\spec R'$, 
where $R'=k[[x_1,\ldots, x_N]]/(I+(g_1,\ldots, g_c))$. 
Under the notation in \ref{notation}, we obtain
$$
{\cal O}_{(X')_m^0} = {\cal O}_{(X)_m^0} \ / \ ( \{(G_{i})_1^0,
\ldots , (G_{i})_m^0\}_{i=1}^c ),
$$
where we identify $(G_{i})_n^0 \in k[{\underline {X}}_1, \ldots ,
{\underline {X}}_m]$ with its class in ${\cal O}_{(X)_m^0}$. Since
${\cal O}_{X_m^0}$ is a catenary ring, applying Krull's theorem,
we obtain
\begin{equation}
\label{15}
\dim {\cal O}_{X_m^0}=\dim {\cal O}_{(X')_m^0} + \text{ht } (
\{(G_{i})_1^0, \ldots , (G_{i})_m^0\}_{i=1}^c ) \leq \dim {\cal
O}_{(X')_m^0}  + m c. 
\end{equation}

Therefore, from (\ref{formula}) we have $$\widehat{\mld}(x_0;X, {\cal J}_X)
=\inf_m\{(m+1)d-\dim(\spec {\widehat R})^0_m\}$$
$$ \geq
\inf_m\{(m+1)(d-c)-\dim(\spec { R'})^0_m\}+c
= (d-c-1)+c =d-1.$$
In addition, $X$ has a
singularity at $x_0$, because $(X', x_0)$ is singular, hence
$\widehat{\mld}(x_0;X, {\cal
J}_X) =d-1$ by Proposition \ref{basic}. $\Box$




\begin{Lemma} 
\label{hyper} If $X$ has a top singularity at $x_0$, then $X$ is
locally  a hypersurface of multiplicity $2$ at $x_0$.
\end{Lemma}

{\it Proof:} From (\ref{formula}) it follows that a germ of an algebraic variety
$(X,x_0)$ is a top singularity if and only if
\begin{equation}
\label{16}
\dim X_m^0 \leq m d + 1 \ \ \ \ \ \ \text{ for every } m \geq 1
\end{equation}

and the equality holds at least for an integer $m$. Suppose that
$(X,x_0)$ is a top singularity. We may suppose that $X$ is affine.
Let $N$ be the embedding dimension of $X$ at $x_0$, then $N \geq
d+1$ because $x_0$ is a singular point of $X$. Besides, with the
notation in 2.2, we have $X_1^0= \text{Spec } k[{\underline
{X}}_1]$; hence, $\dim X_1^0 =N$. Thus,  inequality (\ref{16}) for
$m=1$ implies that $N=d+1$, {\it i.e.,} $X$ is locally  a
hypersurface at $x_0$. Let $X$ be defined by $f(x_1, \ldots, x_{d+1})=0$,
then,
$$
{X_2^0}= \text{Spec } k[{\underline {X}}_1, {\underline {X}}_2] \
/ \ (F^0_2),
$$
and we have $\dim X_2^0 =2d+1$ (resp. $\dim X_2^0 =2(d+1)$) if
$F^0_2  \neq 0$ (resp. $F^0_2 = 0$). Therefore, inequality
(\ref{16}) for $m=2$ implies that $F^0_2 \neq 0$, {\it i.e.},
${\mult}_{x_0} f =2$. $\Box$

\begin{Corollary}
A singularity $(X,x_0)$ is a top singularity if and only if
\begin{equation}
\label{17}
\dim X_m^0 = m d + 1 \ \ \ \ \ \ \text{ for every } m \geq 1 .
\end{equation}
\end{Corollary}

{\it Proof:} The sufficient condition is clear. For the necessary
one, suppose that $(X,x_0)$ is a top singularity. Then, $X$ is
locally  a hypersuperface of multiplicity $2$ at $x_0$, and let it be
defined by $f(x_1, \ldots, x_{d+1})=0$. Then, $\dim X_{1}^0=d+1$
and for  $m \geq 2$, we have
$$
{\cal O}_{X_{m}^0}={\cal O}_{X_{m-1}^0}[{\underline X}_{m}] \ / \
(F^0_{m}),
$$
and hence, $\dim X_{m}^0= \dim X_{m-1}^0 +d + \delta_{m}$, where
$\delta_{m}=0$ (resp. $\delta_{m}=1$) if $F^0_{m}$ is not a zero
divisor (resp. $F^0_{m}$ is a zero divisor) in ${\cal
O}_{X_{m-1}^0}[{\underline X}_{m}]$. Therefore, $\dim X_{m}^0= m d
+ 1 + \sum_{r=2}^m \delta_{r}$. But $(X,x_0)$, a top singularity,
implies that (\ref{16}) holds, and hence $\sum_{r=2}^m \delta_{r} \leq 0$.
Thus,
$\sum_{r=2}^m \delta_{r} = 0$ and (\ref{17}) holds. $\Box$

\begin{defn}
Let $X$ be the hypersurface defined by $x_1 x_2=0$ in $\AA^{d+1} (d\geq 1)$, where $\{x_1, x_{2}\} $ is a part of  the coordinate system of $\AA^{d+1}$.
Then, the singularity $(X, 0)$ is called a {\it normal crossing double}  singularity (sometimes we call it an ncd singularity).

Let $X$ be the hypersurface defined by $x_1^2-x_2^2x_3=0$ in $\AA^{d+1} (d\geq 1)$.
Then, the singularity $(X,0)$ is called a {\it pinch point}.
\end{defn}

\begin{exmp}
\label{ncd}
 Next we give an example of a top singularity of dimension
$d=1$: Let $X$ be a plane curve with an ordinary node, {\it i.e.,} locally it is defined  by $x_1 x_2=0$ in
$\AA^2_k$, and let us consider its germ $(X,{\underline 0})$ at
$\underline 0$.
Then, by inversion of adjunction, we have
$$\hmld(\underline{0}; X, \j_X)=\mld(\underline{0}; \AA^2, I_X)=\mld(\underline{0}; \AA^2, X).$$
It is well known that the right hand side is 0, {\it i.e., } $(X,\underline{0})$ is a top singularity.

We give here another proof by using jet schemes.
 For $m \geq 0$, we have
$$
{X_m^0}= \text{Spec } k[{\underline {X}}_1, \ldots , {\underline
{X}}_m] \ / \ (\{ \sum_{1 \leq i \leq n-1} X_{1,i} X_{2,n-i}\}_{1
\leq n \leq m}).
$$
It follows that $X_m^0$ has $m$ irreducible components given by
$$
X_{1,1}=X_{1,2 }= \ldots =X_{1,r_1}=X_{2,1}= \ldots =X_{2,r_2}=0 \
\ \  \text{for } r_1,r_2 \geq 0, \ r_1+r_2=m-1.
$$
Thus, each irreducible component has dimension $2m-(m-1)=m+1$, and
hence, $\dim X_m^0 = m+1$ for $m \geq 0$. Therefore, $(X,
{\underline 0})$ is a top singularity.
\end{exmp}

\begin{exmp} 
\label{pinch}
 We present a two-dimensional example in the following:
Let $X\subset \AA^3$ be the hypersurface defined by $x_1^2-x_2^2x_3=0$.
Then, $(X,\underline{0})$ is a top singularity.
Indeed, let $\varphi: A'\to A=\AA^3$ be the blow-up at the singular locus 
 of $X$ and let $E$ be the exceptional divisor for $\varphi$, then the strict transform $Y$ of $X$ in $A'$ is non-singular and crosses $E$ normally.
By the inversion of adjunction, we have
$$\hmld(\underline{0}; X, \j_X)=\mld(\underline{0};A, I_X).$$
Here as $A'$ is also non-singular, the right hand side is 
$$\mld(\underline{0};A, I_X)=\mld(\underline{0};A, X)=\mld (\varphi^{-1}(\underline{0}); A', Y+E),$$
because $I_X\o_{A'}=\o_A'(-Y-2E)$, $K_{A'/A}=E$ and $K_{{\overline A}/A}=K_{{\overline A}/A'}+f^*(K_{A'/A})$, where $f:\overline A \to A'$ is a resolution.
As $\varphi^{-1}(\underline{0})$ is a curve that does not contain the double locus of $Y+E$,
it is well known that 
$$\mld (\varphi^{-1}(\underline{0}); A', Y+E)=1.$$

\end{exmp}

\begin{Proposition}
\label{product}
  A normal crossing double singularity and a pinch point are top singularities. 
\end{Proposition}

{\it Proof:} From 
Example \ref{ncd}. for $d=1$, an ncd singularity $(X, 0)=(C, 0)$ is a top singularity.
Let $d\geq 2$ and $(X,0)$ be a $d$-dimensional ncd singularity.
Then, the $(d-1)$ successive general hyperplane cut gives an ordinary double point $(C, 0)$
of a curve $C$.
Then, by Lemma \ref{ci},  $(X, 0)$ is a top singularity. 

For the pinch point of dimension $d$, by the same argument as above, we can reduce 
the discussion to the  two-dimensional case  (Example \ref{pinch}).
 $\Box$

\begin{rem} For $d \geq 2$, recall that if $X
\subset \AA^{d+1}_k$ is a normal hypersurface and $x_0$ is a closed
point of $X$, then
$$
\aligned \widehat{\mld}(x_0;X, {\cal J}_X)  = {\mld}(x_0;X,
{\cal O}_X) = \ \ \ \ \ \ \ \ \ \ \ \ \ \ \ \ \ \ \ \ \ \ \ \ \ \
\ \ \ \ \ \ \ \ \ \ \ \ \ \ \ \ \ \ \ \ \\
\ \ \ \ \ \ \ \ \ \ \ =\inf \{k_E +1 \mid  \  \nu_E \text{
divisorial valuation centered at } x_0 \}. 
\endaligned
$$
Note that if $\pi:Y \rightarrow X$ and $\pi':Y' \rightarrow  X$
are two desingularizations of $X$ and $Y'$ dominates $Y$, let
$\rho:Y' \rightarrow Y$ be such that $\pi'=\rho \circ \pi$, then
we have $ K_{Y'/X}= K_{Y'/Y} + \rho^*(K_{Y/X})$ and $K_{Y'/Y} $ is
effective. Therefore, given a normal hypersurface $X \subset
\AA^{d+1}_k$ of dimension  $d \geq 2$,  to prove that
$(X,x_0)$ is a top singularity, it suffices to show that there
exists a desingularization $\pi:Y \rightarrow X$ such that
\begin{equation}
\label{18}
\inf \{k_E +1 \mid  \  E \text{ prime divisor on } Y \text{ such
that } \pi(E)= x_0 \} = d-1 . 
\end{equation}
\end{rem}

\begin{exmp}
\label{RDP}
 Equality (\ref{18}) is satisfied for the minimal desingularizations of all rational double
points of dimension $2$ (also called Du Val singularities) because they are
canonical singularities of dimension $d=2$. The following is a
list of rational double points, for each of them, the completion
$\widehat{{\cal O}_{X,{\underline 0}}}$ of the local ring ${\cal
O}_{X,{\underline 0}}$ of its germ at $\underline 0$ is described
as a quotient of the ring of series $k[[x_1, x_2, x_3]]$. More
precisely, for each of the types of the rational double points on
the left hand side, there exist ${\bf x}_1, {\bf x}_2, {\bf x}_3$
generating the maximal ideal of $\widehat{{\cal O}_{X,{\underline
0}}}$ and satisfying the equation on the right hand side (recall
that $\text{char } k =0$):
$$
\aligned
{\bf A}_n (n \geq 1): \ \ \ \ \ \  {\bf x}_1^2 + {\bf x}_2^2 + {\bf x}_3^{n+1} =0 \ \ \ \\
\ \ \ \ \ \ {\bf D}_n (n \geq 4): \ \ \ \ \ \ {\bf x}_1^2 + {\bf x}_2^2 {\bf x}_3 + {\bf x}_3^{n-1} =0 \\
{\bf E}_6 : \ \ \ \ \ \ \ \ \ \ \ \ \ \ \ \ \ {\bf x}_1^2 + {\bf x}_2^3 + {\bf x}_3^4 =0 \ \ \ \ \ \ \\
{\bf E}_7 : \ \ \ \ \ \ \ \ \ \ \ \ \ \ \ \ \ {\bf x}_1^2 + {\bf x}_2^3 + {\bf x}_2 {\bf x}_3^3 =0 \ \ \ \\
{\bf E}_8 : \ \ \ \ \ \ \ \ \ \ \ \ \ \ \ \ \  {\bf x}_1^2 + {\bf
x}_2^3 + {\bf x}_3^5 =0\ \ \  \ \ \
\endaligned
$$
\end{exmp}

\begin{Proposition}
\label{cDV-top}
 A compound Du Val singularity is a top singularity.
\end{Proposition}

{\it Proof:} 
Let $(X,x_0)$ be a compound Du Val singularity of dimension $d\geq 3$.
Then, a successive $(d-2)$ hyperplane cut produces a Du Val singularity.
As in the previous example, Du Val singularities are top singularities.
By Lemma \ref{ci}, we obtain that $(X,x_0)$ is a top singularity. $\Box$

\begin{say}
We will prove that  cDV, ncd, and pinch points 
are all top singularities.
Recall that given
$f(x_1, \ldots , x_{d+1}) \in k[x_1, \ldots , x_{d+1}]$ (resp. 
$f
\in k[[x_1, \ldots , x_{d+1}]]$), if $\text{in} f$ denotes the
initial form of $f$ in the graded ring $k[x_1, \ldots , x_{d+1}]$
(resp. $k[[x_1, \ldots, x_{d+1}]]$), with the usual
graduation, then the smallest possible dimension $\tau$
of a linear subspace $V_0$ of $V=k x_1+ \ldots + k x_{d+1}$ such
that $\text{in} f $ lies in the subalgebra $k[V_0]$ of $k[x_1,
\ldots , x_{d+1}]$ is an invariant of the germ $(X, {\underline
0})$ of the hypersurface $X \subset \AA^{d+1}_k$  at $\underline 0$
defined by $f(x_1, \ldots , x_{d+1})=0$ (resp. of $\text{Spec }
k[[x_1, \ldots , x_{d+1}]] / (f)$  ) ([Hi1], chap. III). We denote
it by $\tau(X, {\underline 0})$ (resp. by $\tau(f)$). Given a germ
$(X,x_0)$ of a hypersurface in $\AA^{d+1}_k$  at a closed point
$x_0$, the $\tau$-invariant $\tau(X, x_0)$ is defined as the
$\tau$-invariant of the germ of a hypersurface obtained after a
translation of $x_0$ to $\underline
0$. \\
\end{say}

\begin{Lemma}
\label{tau}
Let $(X, {\underline 0})$ be the germ  of a
hypersurface $X \subset \AA^{d+1}_k$ of multiplicity $2$ at $\underline 0$. Then,
there exist ${\bf x}_1, \ldots, {\bf x}_{d+1} \in \widehat{{\cal
O}_{X, {\underline 0}}}$ generating its maximal ideal and satisfying
$$
{\bf x}_1^2+ \ldots + {\bf x}_{\tau}^2 + g({\bf x}_{\tau+1},
\ldots , {\bf x}_{d+1}) =0
$$
with $\tau=\tau(X, {\underline 0})$, $g(x_{\tau+1}, \ldots ,
x_{d+1}) \in k[[x_{\tau+1}, \ldots , x_{d+1}]]$ and either  $g=0$
or ${\mult} \ g \geq 3$.
\end{Lemma}

{\it Proof:} 
This is well known (for example, see [KM] 4.24 and 4.25).
Actually, iterating the procedures Steps 1 --  3 in 4.25 of [KM],
we obtain the required equation.
 $\Box$
 
\begin{Proposition}
\label{tau>1}
Let $(X, {x_0})$ be a germ of a hypersurface of multiplicity $2$
and $\tau (X, {x_0}) >1$. Then, $(X,{x_0})$ is either an ncd singularity or a cDV singularity; therefore, it is a top singularity.
\end{Proposition}

{\it Proof:} As in Lemma \ref{tau}, let ${\bf x}_1, \ldots, {\bf x}_{d+1} \in \widehat{{\cal
O}_{X, {\underline 0}}}$  generate its maximal ideal and satisfy
$$
{\bf x}_1^2+ \ldots + {\bf x}_{\tau}^2 + g({\bf x}_{\tau+1},
\ldots , {\bf x}_{d+1}) =0
$$
with $\tau=\tau(X, {\underline 0})$, $g(x_{\tau+1}, \ldots ,
x_{d+1}) \in k[[x_{\tau+1}, \ldots , x_{d+1}]]$ and either $g=0$
or ${\mult} \ g \geq 3$.

 If
$\tau \geq 3$, let $R':={\widehat R}/({\bf x}_4,\ldots, {\bf x}_{d+1})$,
where ${\bf x}_i$ denotes the class of $x_i$ in $\widehat R=\widehat{\o_{X,\underline 0}}$.
Then, $\spec R'$ is defined in $\spec k[[x_1,x_2,x_3]]$ by
$$
x_1^2 \ + \ x_2^2 \ + \ x_3^2=0.
$$
Hence, it has an ${\bf A}_1$-singularity at ${\underline 0}$,  and therefore, 
$(X,x_0)$ is a cDV singularity (Corollary \ref{3.4}).  

If $\tau=2$ and $g=0$, then $(X, x_0)$ is an ncd singularity.

If $\tau=2$ and $g\neq 0$, then 
 there exists ${\underline
\lambda}=(\lambda_4, \ldots, \lambda_{d+1}) \in \AA^{d-1}_k$ such
that $g(x_3, \lambda_4 x_3, \ldots,  \lambda_{d+1} x_3)$ is
non-zero and  its multiplicity is
$m={\mult}_{\underline 0} \ g(x_3, \ldots , x_{d+1})$. 
Let $R':={\widehat R}/({\bf x}_4-\lambda_4{\bf x}_3,\ldots, {\bf x}_{d+1}-\lambda_{d+1}{\bf x}_{3})$, then $\spec R'$ is defined in $\spec k[[x_1,x_2,x_3]]$ by
$$
x_1^2 \ + \ x_2^2 \ + \ u \ x_3^m=0,
$$
where $u$ is a unit in $k[[x_3]]$; hence, $\spec R'$ has
an ${\bf A}_{m-1}$-singularity (see Example \ref{RDP}); thus, $(X, x_0)$ is a
cDV 
singularity. $\Box$

\begin{say} 
\label{tau=1} Let $(X, {\underline 0})$ be a germ of a hypersurface $X
\subseteq \AA^{d+1}_k$ of multiplicity $2$ and $\tau(X,
{\underline 0})=1$. Let ${\bf x}_1, \ldots, {\bf x}_{d+1}$ be
generating the maximal ideal of $\widehat{{\cal O}_{X, {\underline
0}}}$ and satisfying
\begin{equation}
\label{20}
{\bf x}_1^2+ g({\bf x}_{2}, \ldots , {\bf x}_{d+1}) =0 ,
\end{equation}
where $g(x_{2}, \ldots , x_{d+1}) \in k[[x_{2}, \ldots ,
x_{d+1}]]$ and since $X$ is reduced, $g \neq 0$ and  ${\mult}
\ g \geq 3$ (Lemma \ref{tau}). Let us consider the germ  of the hypersurface \linebreak $g(x_{2}, \ldots , x_{d+1})=0$ at ${\underline
0}$
in $\text{Spec } k[[x_1, \ldots, x_{d+1}]]$. Although this germ
depends on the choice of  $x_1, \ldots, x_{d+1}$, its multiplicity
$m_2:={\mult} \ g$, and its $\tau$-invariant at ${\underline
0}$, let it be $\tau_2$, which only depends on $(X,{\underline 0})$ (this
follows from [Hi2]. See Remark \ref{hiro} ). Given a germ $(X,x_0)$
of a hypersurface in $\AA^{d+1}_k$ at a closed point $x_0$, we
define $m_2(X,x_0)$ and $\tau_2(X, x_0)$ to be the invariants
defined as before, after a translation of $x_0$
to $\underline 0$.
\end{say}

\begin{rem}
 \label{hiro}
  In [Hi2], the following combinatorial object has been
 considered: Given $f \in k[[x_1, \ldots, x_{d+1}]]$, let
 $f=\sum_{n_1, {\underline n}} c_{n_1, {\underline n}} \ x_1^{n_1}
 {\underline x}^{\underline n}$, where $c_{n_1, {\underline n}} \in
 k$, \linebreak ${\underline x}=(x_2, \ldots, x_{d+1})$ and $(n_1,
 {\underline n})$ runs in $\ZZ_{\geq 0} \times (\ZZ_{\geq 0})^d$,
 let $\text{Supp} (f; {\underline x};x_1) = \{(n_1, {\underline
 n})$ $\text{such that } c_{n_1, {\underline n}} \neq 0 \}$, and
 let $p_m: \RR \times \RR^{d} \rightarrow \RR^d$ be the projection
 from $(m,{\underline 0})$, where $m:=\text{mult} f$. Then
 $$
 p_m(\text{Supp} (f; {\underline x};x_1)) \ = \ \left\{
 \frac{1}{m-n_1} \ {\underline n} \ / \ c_{n_1, {\underline n}}
 \neq 0 \right\}.
 $$
 Let $\Delta(f; {\underline x};x_1)$ be the convex hull of
 $p_m(\text{Supp} (f; {\underline x};x_1)) + (\RR_{\geq 0})^d$, and
 let
 $$
 \Delta(f; {\underline x}):= \cap_{x'} \Delta(f; {\underline x};x'),
 $$
 where $x'$ runs in the set of elements of $k[[x_1, \ldots,
 x_{d+1}]]$ such that $x', x_2, \ldots x_{d+1}$ is a regular system
 of parameters of $k[[x_1, \ldots, x_{d+1}]]$. The combinatorial
 object $\Delta(f; {\underline x})$ is called the {\it (first)
 characteristic polyhedron of $f$ with respect to $\underline x$}
 ([Hi2], def. 1.2).

 Now, let $f=x_1^2 + g(x_2, \ldots, x_{d+1})$, where $g \in k[[x_2,
 \ldots, x_{d+1}]]$, $g \neq 0$,  and $\text{mult} \ g \geq 3$ (see
 (20)). Because there is no term in $x_1$ in the previous expression
 of $f$, we have $\Delta(f; {\underline x};x_1)=\Delta(f;
 {\underline x})$ ([Hi2] Theorem 4.8). Suppose that $x'_1, \ldots,
 x'_{d+1}$ is another regular system of parameters of $k[[x_1,
 \ldots, x_{d+1}]]$ such that
\begin{equation}
\label{21}
 x_1^2 + g(x_2, \ldots, x_{d+1})= u \ ( (x'_1)^2 + g(x'_2, \ldots,
 x'_{d+1})),  
 \end{equation}
 where $u$ is a unit in $k[[x_1, \ldots, x_{d+1}]]$ and $g' \in
 k[[x'_2, \ldots, x'_{d+1}]]$, $\text{mult} \ g' \geq 3$. Then,
 considering the initial forms in the graded ring $k[[x_1, \ldots,
 x_{d+1}]]$, with the usual graduation, it follows that $x'_1= v
 (x_1 + h)$, where $v$ is a unit and $h \in (x_1, \ldots,
 x_{d+1})^2$. On the other hand, by the definition of $\Delta(f;
 {\underline x};x_1)$ one can naturally consider the graduation on
 $k[[x_1, \ldots, x_{d+1}]]$
 defined by the monomial valuation on $k[[x_1, \ldots, x_{d+1}]]$
 given  by $\nu(x_1)=\frac{1}{2}, \ \nu(x_i)=\frac{1}{\text{mult} \
 g}$ for $2 \leq i \leq d+1$. Then, considering the initial form
 with respect to this graduation on both sides of (\ref{21}), we obtain
 that $\nu(h)> \frac{1}{2}$, $\nu(x'_1)=\nu(x_1)$, and
 $\text{in}_\nu g= \text{in}_\nu g'$, which implies $\text{mult} \
 g= \text{mult} \ g'$ and
 $\tau(g)=\tau(g')$.\
  \end{rem}

\begin{Lemma}
\label{t1m3}
Let $(X, {\underline 0})$ be a germ of a hypersurface in
$\AA^{d+1}_k$ of multiplicity $2$ and $\tau(X, {\underline 0})
=1$. If  $(X, {\underline 0})$ is a top singularity then $d \geq
2$ and $m_2(X,{\underline 0})=3$.
\end{Lemma}

{\it Proof:} Let ${\bf x}_1, \ldots, {\bf x}_{d+1}$ be generating
the maximal ideal of  $\widehat{{\cal O}_{X, {\underline 0}}}$ and
$g \in k[[x_{2}, \ldots , x_{d+1}]]$ such that (\ref{20}) holds. Then,
we have
$$
\aligned X_1^0 \cong \text{Spec } k[{\underline {X}}_1] \ \ \ \ \
\ {X_2^0} \cong \text{Spec } k[{\underline {X}}_1, {\underline
{X}}_2] \
/ \ ((X_{1,1})^2)  \\
{X_3^0} \cong \text{Spec } k[{\underline {X}}_1, {\underline
{X}}_2, {\underline {X}}_3] \ / \ ((X_{1,1})^2, G_3^0) \ \ \ \ \ \
\ \ \
\endaligned
$$
(see (\ref{3}) and (\ref{4})). Moreover, $X_{1,1}$ does not divide $G_3^0$;
hence, $\dim X_3^0 =3d+1$ (resp. $\dim X_3^0 =3d+2$) if $G^0_3 \neq
0$ (resp. $G^0_3 = 0$). Therefore, if $(X, {\underline 0})$ is a
top singularity then $G^0_3 \neq 0$ (see (\ref{17})) or
equivalently,
$m_2(X,{\underline 0})={\mult}_{\underline 0} g=3$.

Now let us show that if $(X, {\underline 0})$ is a top
singularity, then $d \geq 2$. In fact, suppose that $d=1$ and
$m_2(X,{\underline 0})=3$. Then, there exist ${\bf x}_1, {\bf x}_2$
generating the maximal ideal of $\widehat{{\cal O}_{X, {\underline
0}}}$  such that ${\bf x}_1^2+ {\bf x}_2^3=0$. In fact, there
exist ${\bf x}_1, {\bf x}_2 \in \widehat{{\cal O}_{X, {\underline
0}}}$ such that ${\bf x}_1^2 + g({\bf x}_2)=0$, where $g(x_2) \in
k[[x_2]]$ has multiplicity $3$; hence, $g(x_2)=x_2^3 \ u(x_2)$,
where $u(x_2)$ is a unit in $k[[x_2]]$. Because $k$ is algebraically
closed, there exists $v(x_2) \in k[[x_2]]$ such that
$v(x_2)^3=u(x_2)$. Then, replacing ${\bf x}_2$ by $v({\bf x}_2) \
{\bf x}_2$, we have ${\bf x}_1^2+ {\bf x}_2^3=0$. Therefore,
$$
\aligned ({X_4^0})_\text{red} \cong \text{Spec } k[{\underline
{X}}_1, \ldots, {\underline {X}}_4] \ / \ (X_{1,1}, X_{2,1},
X_{1,2})
\\
({X_5^0})_\text{red} \cong \text{Spec } k[{\underline {X}}_1,
\ldots, {\underline {X}}_5] \ / \ (X_{1,1}, X_{2,1}, X_{1,2}).
\endaligned
$$
Hence, $\dim {X_5^0} =7 >5+1$ and $(X, {\underline 0})$ is not
a top singularity. This concludes the proof. $\Box$

\begin{Proposition}
\label{t1m3t2}
Let $(X, {x_0})$ be a germ of a  hypersurface in $\AA^{d+1}_k$ of
multiplicity $2$ and $\tau (X, {x_0}) =1$. If $m_2(X,{x_0})=3$ and
$\tau_2(X,{x_0}) >1$, then $(X,{x_0})$ is either a cDV singularity or a pinch point; therefore, it is a top singularity.
\end{Proposition}

{\it Proof:} We may suppose that the point $x_0$ is the origin
${\underline 0} \in \AA^{d+1}_k$. 
 We can take generators ${\bf x}_1, \ldots, {\bf x}_{d+1}$ 
 of  the maximal ideal of $\widehat{{\cal O}_{X, {\underline
0}}}$  such that
$$
{\bf x}_1^2+ g_3({\bf x}_{2},\ldots, {\bf x}_{\tau_2+1})+h({\bf x}_2\ldots , {\bf x}_{d+1}) =0 .
$$
Here $g_3(x_{2}, \ldots , x_{\tau_2+1}) \in k[x_{2}, \ldots , x_{\tau_2+1}]$ is homogeneous  of degree $3$
and $h\in k[[x_{2}, \ldots,
x_{d+1}]]$ has multiplicity $\geq 4$.
Let $R'=\widehat R/({\bf x}_{\tau_2+2},\ldots, {\bf x}_{d+1})$, 
then $\spec R'$ is defined in $\spec k[[x_1,\ldots, x_{\tau_2+1}]] $ by

\begin{equation}
\label{tau2}
{ x}_1^2+ g_3({ x}_{2},\ldots, { x}_{\tau_2+1})+g_4({x}_2\ldots , { x}_{\tau_2+1}) =0 ,
\end{equation}
where $g_3$ is homogeneous of degree 3 as stated before and $\mult g_4\geq 4$.
Next, take a general combination $\lambda=(\lambda_2,\ldots, \lambda_{\tau_2+1})\in \AA_k^{\tau_2}$
so that the hyperplane $H$ in $\PP_k^{\tau_2-1}$ defined by $\sum_{i=2}^{\tau_2+1}\lambda_ix_i=0$ does not contain any irreducible components of
 the hypersurface $G$ defined by $g_3=0$. 
Here we  note that the hypersurface $G$ is not a triple plane because $\tau_2>1$ and 
 by the definition of $\tau_2$. 
Therefore, we may also assume that $G\cap H$ is not a triple hyperplane in $H$.
Now let $R''=R'/(\sum_{i=2}^{\tau_2+1}\lambda_i{\bf x}_i)$.
As $\lambda$ is general, we may assume that $\lambda_{\tau_2+1}\neq 0$, and therefore, we can eliminate $x_{\tau_2+1}$, {\it i.e.,} $\spec R''$ is defined in 
$\spec k[[x_1,\ldots, x_{\tau_2}]]$ by
$$
{ x}_1^2+ g'_3({ x}_{2},\ldots, { x}_{\tau_2})+g'_4({ x}_2\ldots , { x}_{\tau_2}) =0 ,$$
where $g'_3$ and $g'_4$ have the same properties as $g_3$ and $g_4$ in (\ref{tau2}),
respectively.
By reiterating this argument, we obtain a two-dimensional singularity whose germ is 
\begin{equation}
\label{26}
\widehat{R} \ \cong \ k[[x_1,x_2,x_3]] \ / \ ( \ x_1^2 \ + \
g(x_2,x_3) \ ),
\end{equation}
where $g(x_2,x_3) \in k [[x_2,x_3]]$, ${\mult} \ g(x_2,x_3)=3$,
and by the above argument,
$\text{in} g(x_2,x_3)$ does not give a triple point in $\PP_k^1$, {\it i.e.},
$\text{in} g(x_2,x_3)$ has at
least two different linear factors. 
In this case, after a possible change of the regular system of parameters of 
$k[[x_2,x_3]]$, we may suppose that $g(x_2,x_3)=x_3(x_2^2 + v x_3^m)$,  
where $m \geq 2$ and $v \in k[[x_3]]$ is either a unit or $0$ ([KM], 
Step 4 in 4.25).
Here if $v$ is a unit, then $\spec \widehat R $ has ${\bf D}_{m+2}$-singularity and if $v=0$,
then $\spec \widehat R $ has a pinch point.
$\Box$

\begin{say}
\label{3.9}
Let $(X, {\underline 0})$ be a germ of a hypersurface $X
\subseteq \AA^{d+1}_k$ of multiplicity $2$, $\tau(X, {\underline
0})=1$, $m_2(X,{\underline 0})=3$, and $\tau_2(X,{\underline
0})=1$. Then, there exist ${\bf x}_1, \ldots, {\bf x}_{d+1}$
generating the maximal ideal of $\widehat{{\cal O}_{X, {\underline
0}}}$  such that
\begin{equation}
\label{27}
{\bf x}_1^2 \ + \ {\bf x}_2^3 \ +  \ g_3({\bf x}_3, \ldots, {\bf
x}_{d+1}) \ {\bf x}_2 \ + \ g_4 ({\bf x}_3, \ldots, {\bf x}_{d+1})
=0, 
\end{equation}
where $g_i \in k[[x_3, \ldots, x_{d+1}]]$ and
${\mult}_{\underline 0} \ g_i \geq i$, for $i=3,4$. In fact,
there exist $x_1, \ldots, x_{d+1}$ whose classes in
$\widehat{{\cal O}_{X, {\underline 0}}}$ generate the maximal
ideal, and $g \in k[[x_2, \ldots , x_{d+1}]]$ such that (\ref{20})
holds.   Moreover, because ${\mult} \ g=m_2(X, {x_0})=3$ and
$\tau_2(X, {x_0})=1$, by Weierstrass' preparation theorem and
after a Tschirnhausen transformation, we may suppose that
$$
g(x_2, \ldots , x_{d+1})= u \ \left( x_2^3 \ + \ g_3 (x_3, \ldots,
x_{d+1}) \ x_2 \ + \ g_4(x_3, \ldots, x_{d+1}) \right),
$$
where $u$ is a unit in $k[[x_2, \ldots , x_{d+1}]]$ and $g_i
\in k[[x_3, \ldots , x_{d+1}]]$ is such that ${\mult} \ g_i
\geq i$ for $i=3, 4$. Replacing $x_1$ by $v x_1$, where $v$ is a
unit in $k[[x_2, \ldots , x_{d+1}]]$ such that $v^2=u$, and
considering the equality induced on the classes ${\bf x}_i$ of
$x_i$ in $\widehat{{\cal
O}_{X,{\underline 0}}}$,  we obtain (\ref{27}).

Given $g_3, g_4 \in k[[x_3, \ldots, x_{d+1}]]$  as in (\ref{27}), let
$$
m_3(g_3, g_4):= 6 \ \text{min} \left\{ \frac{\text{mult} \
g_3}{2}, \ \frac{\text{mult} \ g_4}{3} \right\}.
$$
Note that $m_3(g_3, g_4)\in \NN \cup \{\infty\}$ and $\frac{1}{6}
\ m_3(g_3, g_4) >1$. Moreover, $m_3(g_3, g_4)$ is an invariant of
$(X,{\underline 0})$. 
This follows from [Hi2] (see Remark \ref{hiro}).
Let  $m_3(X,{\underline 0})$ denote this invariant.

\end{say}

\begin{Proposition}
\label{t1m3t1}
Let $(X, {x_0})$ be a germ of a  hypersurface in $\AA^{d+1}_k$ of
multiplicity $2$ and $\tau(X, {x_0}) =1$, $m_2(X,{x_0})=3$, and
$\tau_2(X,{x_0}) =1$. Then, 
the following are equivalent:
\begin{enumerate}
\item[(i)]
 $(X,{x_0})$ is a top singularity,
 \item[(ii)]
there exist ${\bf x}_1, \ldots, {\bf x}_{d+1}$ generating the
maximal ideal of $\widehat{{\cal O}_{X,x_0}}$ such that
$$
{\bf x}_1^2 \ + \ {\bf x}_2^3 + g_3({\bf x}_3, \ldots, {\bf
x}_{d+1}) \ {\bf x}_2 \ + \ g_4({\bf x}_3, \ldots, {\bf x}_{d+1})
\ = \ 0.
$$
where $g_i \in k[[x_3, \ldots, x_{d+1}]]$,  ${\mult} \ g_i
\geq i$ for $i=3,4$ and either ${\mult} \ g_3=3$ or $4 \leq
{\mult} \ g_4 \leq 5$,
\item[(iii)] $\frac{1}{6}m_3(X, x_0)<2$,
\item[(iv)]  $(X,{x_0})$ is a cDV singularity.
\end{enumerate}
\end{Proposition}

{\it Proof:} Implication (iv) $\Rightarrow$ (i) is obvious by Proposition \ref{cDV-top}.
We will show (ii) $\Leftrightarrow$ (iii), and then (i) $\Rightarrow$ (ii) and (ii) $\Rightarrow$ (iv).

Let $R={\cal O}_{X,x_0}$, then we have
\begin{equation}
\label{32}
\widehat{R} \ \cong \  k[[x_1, \ldots, x_{d+1}]] \ / \ ( x_1^2 \ +
\ x_2^3 \ + \  g_3 \ x_2 \ + \ g_4 ), 
\end{equation}
where $g_3, g_4 \in k[[x_3, \ldots, x_{d+1}]]$ and
$\mu_i:={\mult} \ g_i \geq i$, for $i=3,4$ (see \ref{3.9}). Note
that $\mu_i=\infty$ iff $g_i=0$.
 Moreover, by the definition of $m_3(X,x_0)$ (see the
paragraph before Proposition \ref{t1m3t1} and recall that $\frac{1}{6}
m_3(X,x_0)>1$), the condition $\frac{1}{6} m_3(X,x_0)<2$ is
equivalent to the assertion that either $\mu_3=3$ or $4 \leq \mu_4
\leq 5$. Thus, (ii) is equivalent to (iii).

 To prove (i) $\Rightarrow$ (ii), let us argue by contradiction. Suppose that $
\mu_3 \geq 4$ and $
\mu_4 \geq 6$. Then, applying  (\ref{3}) and (\ref{4}), we obtain
$$
\aligned X_1^0 \cong \text{Spec } k[{\underline {X}}_1] \ \ \ \ \
\ ({X_2^0})_\text{red} \cong \text{Spec } k[{\underline {X}}_1,
{\underline {X}}_2] \
/ \ (X_{1,1})  \\
({X_3^0})_\text{red} \cong \text{Spec } k[{\underline {X}}_1,
{\underline {X}}_2, {\underline {X}}_3] \ / \ (X_{1,1}, X_{2,1}) \
\ \ \ \ \ \ \ \ \ \ \ \ \ \ \ \\
({X_4^0})_\text{red} \cong \text{Spec } k[{\underline {X}}_1,
\ldots , {\underline {X}}_4 ] \ / \
(X_{1,1}, X_{2,1}, X_{1,2}) \ \ \ \ \ \ \ \ \ \\
({X_5^0})_\text{red} \cong \text{Spec } k[{\underline {X}}_1,
\ldots , {\underline {X}}_5 ] \ / \ (X_{1,1}, X_{2,1}, X_{1,2}). \
\ \ \ \ \ \ \ \
\endaligned
$$
That is, with the notation in \ref{notation}, we have $F_5^0 \in (X_{1,1},
X_{2,1}, X_{1,2})$, where $f=x_1^2  +  x_2^3  +   g_3 \ x_2 + g_4$
(see (\ref{32})). Therefore, $\dim {X_5^0} = 5 d + 2$, and hence,
$(X,{x_0})$ is not a top singularity (see (\ref{16})).\\

To prove (ii) $\Rightarrow$ (iv),
let us suppose that either $\mu_3=3$ or $4 \leq \mu_4 \leq 5$. Let
${\underline \lambda}=(\lambda_4, \ldots, \lambda_{d+1}) \in
\AA^{d-2}_k$ be such that
$$
{\mult} \ g_i(x_3, \lambda_4 x_3, \ldots, \lambda_{d+1} x_3) =
{\mult} \ g_i(x_3,  x_4, \ldots, x_{d+1}) = \mu_i \ \ \
\text{for } i=3,4.
$$
Hence, for $i=3,4$, if $g_i \neq 0$ then
\begin{equation}
\label{33}
g_i(x_3, \lambda_4 x_3, \ldots, \lambda_{d+1} x_3)  = u_i \
x_3^{\mu_i},  
\end{equation}
where $u_i$ is  a unit in $k[[x_3]]$. Let us consider
$$
R':= \widehat{R} \ / \ ({\bf x}_4-\lambda_4 {\bf x}_3,\ \ldots \ ,
{\bf x}_4-\lambda_{d+1} {\bf x}_{d+1}),
$$
where ${\bf x}_i$ is the class of $x_i$ in $\widehat{R}$, $3 \leq
i \leq d+1$.
By Lemma \ref{special}, it is sufficient to prove that $\spec R'$ has a Du Val singularity.


Note that, by (\ref{32}) and (\ref{33}),
$$
R' \ \cong \ k[[x_1, x_2, x_{3}]] \ / \ (x_1^2 \ + \ x_2^3 \ + \
g'_3(x_3) \ x_2 \ + \ g'_4 (x_3)),
$$
where, for $i=3,4$, we have $g'_i(x_3)=0$ if $\mu_i=\infty$ and
$g'_i(x_3)= u_i \ x_3^{\mu_i}$ if $\mu_i < \infty$.
Here $g'_i$ has the same property on the multiplicity as $g_i$ in (ii).
Then, by Steps 5 --  8 of 4.25 in [KM], 
we obtain that $\spec R'$ has ${\bf E}_n$-singularity $(n=6,7,8)$.
 $\Box$

\vskip.5truecm
The following summarizes the discussions of characterization of a top singularity
(Proposition \ref{tau>1}, Lemma \ref{t1m3}, Proposition \ref{t1m3t2}, Proposition \ref{t1m3t1}).

\begin{Corollary}
\label{maincor}
A germ of a variety $(X, x_0)$ of dimension $d=1$ is a top
singularity if and only if it is the germ  of a
hypersurface in $\AA^2_k$ of multiplicity $2$ at $x_0$and
$\tau(X,x_0)>1$.

A germ of a variety $(X, x_0)$ of dimension $d \geq 2$ is a top
singularity if and only if it is the germ  of a
hypersurface in $\AA^{d+1}_k$ of multiplicity $m_1(X,x_0)=2$ at $x_0$ such
that one of the following holds:
\begin{enumerate}
 \item[(i)]
$\tau(X,x_0)>1$
%
$$
(ii)\  \tau(X,x_0)=1, \ \frac{m_2(X,x_0)}{m_1(X,x_0)!}<2 \text{
and } \tau_2(X,x_0)>1 \ \ \ \ \ \ \ \ \ \ \ \ \ \ \ \ \ \ \ \ \ \
\ \ \ \ \ \ \ \
$$
%
$$(iii)\  \tau(X,x_0)=1, \
\frac{m_2(X,x_0)}{m_1(X,x_0)!}<2, \ \tau_2(X,x_0)=1 \text{ and }
\frac{m_3(X,x_0)}{m_2(X,x_0)!}<2. \ \ \ \ \ \ \
$$
\end{enumerate}
\end{Corollary}

{\it Proof:} For $d=1$, the result follows from Proposition \ref{product}
and Lemmas \ref{hyper}  and \ref{t1m3}. For $d \geq 2$, note that, if
$\tau(X,x_0)=1$, then $\frac{m_2(X,x_0)}{m_1(X,x_0)!}= \frac{1}{2}
m_2(X,x_0)$ is well defined and is always $>1$ (see \ref{tau=1} and
Remark \ref{hiro}), hence $\frac{1}{2} m_2(X,x_0) < 2$ is equivalent to
$m_2(X,x_0)=3$. Also note that  if $\tau(X,x_0)=1$, $m_2(X,x_0)=3$
and $\tau_2(X,x_0)=1$, then $\frac{m_3(X,x_0)}{m_2(X,x_0)!}=
\frac{1}{6} m_3(X,x_0)$ is well defined (see \ref{3.9}). Thus, the
result follows from Lemmas  \ref{hyper}, \ref{t1m3}, and Propositions \ref{tau>1}, \ref{t1m3t2},
and \ref{t1m3t1}. $\Box$

\begin{Theorem}
\label{topequation}
A germ of a variety $(X, x_0)$ of dimension $d$ is a top singularity
 if and only if either
\begin{enumerate}
\item[(i)] 
there exist a minimal system of generators ${\bf x}_1, {\bf x}_2,\ldots, {\bf x}_{d+1}$
of the maximal ideal of $\widehat{{\cal O}_{X, x_0}}$ 
such that ${\bf x}_1 {\bf x}_2=0$ or ${\bf x}_1^2- {\bf x}_2^2{\bf x}_3=0$, or
\item[(ii)]
 $d \geq 2$ and there exist a minimal system of  generators ${\bf x}_1, \ldots , {\bf
x}_{d+1}$ of the maximal ideal of $\widehat{{\cal O}_{X,
x_0}}$  such that one of the following holds:
\begin{itemize}
\item[(a)] ${\bf x}_1^2 \ + \ldots + \ {\bf x}_\tau^2 \ + \ g({\bf
x}_{\tau+1}, \ldots, {\bf x}_{d+1})\ = \ 0$ \ \ where $\tau \geq
2$, $ g(x_{\tau+1}, \ldots, x_{d+1}) \in  k[[x_{\tau+1}, \ldots,
x_{d+1}]]$ and
${\mult} \ g \geq 3$. \\
\item[(b)] ${\bf x}_1^2 \ + \ {\bf x}_2^3 \ + \ p({\bf x}_{3},
\ldots, {\bf x}_{d+1}) \ {\bf x}_2 \ + \ q({\bf x}_{3}, \ldots,
{\bf x}_{d+1}) \ = \ 0$ \ \ where $p(x_{3}, \ldots, x_{d+1})$, \
$q(x_{3}, \ldots, x_{d+1}) \in k[[x_{3}, \ldots, x_{d+1}]]$, \
${\mult} \ p \geq 2$, \ ${\mult} \ q \geq 3$, and either
$2 \leq {\mult} \ p \leq 3$ or $3 \leq
{\mult} \ q \leq 5$. \\
\end{itemize}
\end{enumerate}
\end{Theorem}

{\it Proof:} For $d=1$ the statement is clear by Corollary \ref{maincor}. For $d \geq
2$, note that either (i) or (a) in the theorem holds if and only
if (i) in Corollary \ref{maincor} holds and that (b) in the theorem holds if
and only if either (ii) or (iii)
in Corollary \ref{maincor} holds. $\Box$

\begin{Theorem}
\label{top}
A germ of a variety $(X,x_0)$ is a top singularity if and only if either
\begin{enumerate}
\item[(i)]   $(X,x_0)$ is a normal crossing double singularity or a pinch point or
\item[(ii)] $\dim X\geq 2$ and $(X,x_0) $ is a compound Du Val singularity
\end{enumerate}
\end{Theorem}
{\it Proof:}  Conditions (i) and (ii) imply that $(X,x_0)$ is a top singularity by 
Propositions \ref{product} and \ref{cDV-top}.
The converse follows from the fact that a top singularity is a hypersurface double point and, under the classification of the defining equation according the invariants $\tau$ 
and $m$, a class of top singularities  always satisfy 
 condition either (i) or (ii) (Proposition \ref{tau>1}, Lemma \ref{t1m3}, Proposition \ref{t1m3t2}, and Proposition \ref{t1m3t1}). $\Box$



%

\begin{rem} All the results and proofs in this section remain true if
we replace completions by henselizations; that is, if we replace
$\widehat{{\cal O}_{X,x_0}}$ by the henselization ${\cal
O}_{X,x_0}^h$ of the local ring ${\cal O}_{X,x_0}$ and $k[[x_1,
\ldots, x_N]]$ by $k\{x_1, \ldots, x_N\}$ for every $N \geq 1$.
In fact, we have to apply the version of Weierstrass' preparation
theorem for algebraic series and  Hensel's Lemma in Propositions \ref{t1m3t2}
 and  \ref{t1m3t1}. 
 In particular, Theorem \ref{topequation} remains true if we
replace $\widehat{{\cal O}_{X,x_0}}$ by ${\cal O}_{X,x_0}^h$ 
in both (i) and (ii) and we also replace  $k[[x_{\tau+1}, \ldots,
x_{d+1}]]$ (resp. $k[[x_3, \ldots, x_{d+1}]]$) by $k\{x_{\tau+1},
\ldots, x_{d+1}\}$ (resp. by $k\{x_3, \ldots, x_{d+1}\}$) in (a)
(resp. (b)).
\end{rem}


\section{A variant of Shokurov's conjecture}
In this section we consider the Mather version of Shokurov's second conjecture.
Our main result is the following:
\begin{Theorem} A pair $(X, B)$ consisting of an arbitrary variety $X$ and an  effective $\RR$-Cartier divisor $B$ on $X$ satisfies $$\dim X-1\leq {\hmld}(x;X,{{\mathcal J}_X}B) $$ 
if and only if
either 
\begin{enumerate}
\item[(i)]  $B=0$ and $(X,x)$  is a normal crossing double singularity or a pinch point,
\item[(ii)] $B=0$, $\dim X\geq 2$ and $(X,x) $ is a compound Du Val singularity   or
\item[(iii)] $(X,x)$ is non-singular and $0\leq \mult_xB\leq 1$.
\end{enumerate}

\noindent
In  cases $(\mathrm i)$ and $(\mathrm {ii})$, we have ${\hmld}(x;X,{{\mathcal J}_X})= \dim X -1$ and  
 in case $(\mathrm {iii})$, we have 
${\hmld}(x;X,{{\mathcal J}_X}B)=\mld(x; X, B)= \dim X -\mult_xB$ and the minimal  log discrepancy is computed by the exceptional divisor of the first blow-up at $x$.

\end{Theorem}

{\it Proof:} Let $d=\dim X$ and let $(X, B)$ satisfy the condition $d-1\leq {\hmld}(x;X,{{\mathcal J}_X}B) $
at a closed point $x\in X$. 
If $(X, x)$ is singular, then by Proposition \ref{sh1-answer}, we have 
$\hmld(x;X,{{\mathcal J}_X})\leq d-1$ because $\hmld(x;X,{{\mathcal J}_X})$ is an integer
(see \ref{nash}).
If $B\neq 0$ in a neighborhood of $x$, then $${\hmld}(x;X,{{\mathcal J}_X}B)<\hmld(x;X,{{\mathcal J}_X})\leq d-1,$$ in which case $(X, B)$ does not satisfy the condition of the 
theorem. Therefore, if $(X,x)$ is singular, then $B=0$ and $\hmld(x;X,{{\mathcal J}_X})=d-1$, {\it i.e.,} $(X,x)$ is a top singularity.
A top singularity is characterized in Theorem \ref{top} as in (i) and (ii). 

Hence, it is sufficient to characterize a pair $(X, B)$ such that  $X$ is  non-singular and
${\hmld}(x;X,{{\mathcal J}_X}B)={\mld}(x;X,B)\geq d-1$ in terms of  (iii).

If $d=\dim X =1$, then the statement is obvious since ${\mld}(x;X,B)=1-\mult_xB$.

Assume $d=\dim X\geq 2$ and $(X, B)$ satisfies the inequality
${\mld}(x;X,B)\geq d-1$, then  the exceptional divisor $E_1$ of the 
blow-up $\varphi_1: X_1\to X$ of $X$ at $x$ should have the log discrepancy $$k_{E_1}-\ord_{E_1}\varphi_1^* B+1\geq d-1$$
(see \ref{usualmld} and also (\ref{9})).
As $k_E=d-1$ and $\ord_{E_1}\varphi_1^*B=\mult_xB$, this implies $\mult_xB\leq 1$.

Conversely, we assume (iii), that is $\mult_xB\leq 1$.
Under this condition, we check the log discrepancy of every prime divisor over $X$ with the center at $x$.

Let $E$ be a prime divisor over $X$ with the center at $x$,  let $y\in E$ be the generic point, and let $E$ appear in a resolution $f_0:Y\to X$.
Then, by Zariski's result ( for example, see [Ko], VI, 1.3),
we have a sequence of varieties $X_0,X_1, \ldots, X_n$  and birational maps as follows:
\begin{enumerate}
\item[]
$X_0=X$, $f_0=f$.
\item[] If  $f_i: Y\dashrightarrow
 X_i$ is already defined, then let $Z_i\subset X_i$ be the closure of $p_i=f_i(y)$.
Let $X_{i+1}=B_{Z_i}X_i$ and  $f_{i+1}: Y\dashrightarrow
 X_{i+1}$ be the induced map.
\item[] Then, the final birational map $f_{n}:Y\dashrightarrow
 X_n$ is isomorphic at $y$, 
{\it i.e.,} $E$ appears on $X_n$.
\end{enumerate}
Here $B_{Z_i}X_i$ is the blow-up of $X_i$ with the center $Z_i$.
Let $\varphi_i:X_i\to X_{i-1}$ be the blow-up morphism and $E_i\subset X_i$ be the exceptional divisor dominating $Z_i$.
Note that the first blow-up $\varphi_1:X_1\to X_0=X$ is  at the closed point $x$
because the center of $E$ on $X$ is $x$, whereas $f_n$ is isomorphic at the generic points of $E_n$ and $E$.
We also note that $X_i$ and $E_i$ are non-singular at $p_i$ for every $i=1,\ldots n$.
Indeed, this is proved inductively.
As $X_1 $ is the blow-up at a closed point $x=p_0$, $X_1$ and $E_1$ are non-singular at every point.
Suppose $i\geq 2$ and $X_{i-1}$ and $E_{i-1}$ are non-singular at $p_{i-1}$, 
then $X_i$ is the blow-up with the non-singular center when one restricts the morphism on a
neighborhood of $p_{i-1}$. 
As $p_i$ is on the pull back of this neighborhood,  $X_i$ and $E_i$ are non-singular at $p_i$.

Let $B^{(i)}$ be the strict transform of $B$ on $X_i$, then from [Hi1] II sec. 5, Theorem 3
(p.233), we have 

\begin{equation}
\label{26}
\mult_{p_i}B^{(i)}\leq \mult_{p_{i-1}}B^{(i-1)} \ {\mbox{for\ every\ }} i=1,\ldots, n.
\end{equation}

Let $a(E_i, X,B) $ be the discrepancy of $(X, B)$ at the divisor $E_i$, {\it i.e.,} 
$$a(E_i; X, B)=\ord_{E_i}(K_{X_i/X}-\Phi_i^*(B)),$$
where $\Phi_i:X_i\to X$ is the composite $\varphi_1\circ\cdots\circ \varphi_i$.
Note that the log discrepancy of $(X, B)$ at the divisor $E_i$ is $a(E_i; X, B)+1$.

\noindent
{\bf Claim.} For every $i=1,\ldots, n$
$$a(E_i, X,B)\geq 0 \ \ \mbox{and}\ \ a(E_i, X,B)\geq a(E_{i-1}, X,B).$$

By abuse of notation, we denote the strict transform   of $E_{i-1}\subset X_{i-1}$ on 
$X_j$ $(j\geq i)$ by
the same symbol $E_{i-1}$.
Then, we have 
\begin{equation}
\label{nonsing}
\varphi_i^*(E_{i-1})=E_{i-1}+E_i
\end{equation}
by the non-singularity of $X_{i-1}$ and $E_{i-1}$ at $p_{i-1}$ guaranteed in the discussion  above.

Now we prove the claim  by induction on $i$.
First for $i=1$, 
by substituting $K_{X_1/X}=(d-1)E_1$ and $\varphi_1^*(B)=(\mult_xB)E_1+B^{(1)}$ into 
$$a(E_1; X,B)=\ord_{E_1}(K_{X_1/X}-\varphi_1^*(B)),$$ we obtain
$$a(E_1; X,B)=(d-1)-\mult_xB\geq d-2,$$
which is of course non-negative by our assumption $d\geq 2$.

Let $i\geq 2$ and assume that $a(E_j;X,B)\geq 0$ for all $j\leq i-1$ by induction hypothesis.
Then
$$a(E_i; X,B)=\ord_{E_i}(K_{X_i/X}-\Phi_i^*(B))$$
$$\ \ \ \ \ \ \ \ \ \ \ =\ord_{E_i}\left(K_{X_i/X_{i-1}}+\varphi_i^*\left(K_{X_{i-1}/X}-\Phi_{i-1}^*(B)\right)\right)$$
$$\ \ \ \ \ \ \ \ \ \ \ =\ord_{E_i}\left(K_{X_i/X_{i-1}}+\varphi_i^*(\sum_{j\leq i-1}a(E_j, X, B)E_j
-B^{(i-1)})\right)$$
$$\ \ \ \ \ \ \ \ \ \ \ \ \geq\ord_{E_i}(K_{X_i/X_{i-1}}) +a(E_{i-1};X,B) -\mult_{p_{i-1}}B^{(i-1)}.$$
Here we used (\ref{nonsing}) and the hypothesis of the induction.
We may assume that $\codim \overline{\{p_{i-1}\}}\geq 2$, because if 
$\codim \overline{\{p_{i-1}\}}=1$, then $f_{i-1}$ is already isomorphic at the generic point 
$y\in E$. (We may assume that $n$ is taken to be minimal.)
As $\ord_{E_i}(K_{X_i/X_{i-1}})=\codim\overline{\{p_{i-1}\}}-1\geq 1 $ and 
$\mult_{p_{i-1}}B^{(i-1)}\leq\mult_xB\leq 1$ by (\ref{26}), 
we obtain 
$$a(E_i, X,B)\geq a(E_{i-1}, X,B)\geq 0$$ as claimed.

From this,  the log discrepancy at $E$ is $a(E; X, B)+1=a(E_n;X,B)+1\geq a(E_1; X, B)+1=d-\mult_xB\geq d-1$.
Therefore,  the inequality $\dim X-1\leq {\hmld}(x;X,{{\mathcal J}_X}B)=
\mld(x; X,B) $ holds and the minimal log discrepancy $d-\mult_xB$ is computed by the exceptional divisor of the first blow-up at $x$.
$\Box$

As a corollary of the theorem, we have the ``if" part of Conjecture \ref{sho2}:
\begin{Corollary} The inequality $$\dim X-1<\mld (x;X, B)$$ holds
 if  $(X,x)$ is non-singular and 
$mult_x B<1$.
In this case, the minimal log discrepancy is computed by the exceptional divisor of the
first blow-up at $x$.
\end{Corollary}

Since $\mld(x; X, B)=\hmld(x;X,\j_XB)$ for $(X,x)$, which is a normal complete intersection by (\ref{9}),
 we have the following statement for usual mld as a further corollary. 
 
\begin{Corollary} 
A pair $(X, B)$ consisting of a normal complete intersection variety $X$ at a closed point
$x$ and an  effective $\RR$-Cartier divisor $B$ on $X$ satisfies $$\dim X-1\leq {\mld}(x;X,B) $$ 
if and only if
either {\rm (i)} or {\rm(ii)} or {\rm(iii)} in the theorem holds.
\end{Corollary}
By this and Proposition \ref{basic}, Corollary \ref{c.i.} follows.

\vskip.5truecm
\noindent
{\bf References.}
\begin{itemize} \parindent=10mm

\item[\hbox to\parindent{\enskip\mbox{[Am1]}}] F. Ambro, {\em On minimal log discrepancies}, Math. Res. Letters, 6, 573--580 (1999).

\item[\hbox to\parindent{\enskip\mbox{[Am2]}}] F. Ambro, {\em Inversion of adjunction 
for non-degenerate hypersurfaces}, 
Manuscripta Math.  111, 43--49 (2003). 

\item[\hbox to\parindent{\enskip\mbox{[B]}}]
E. Brieskorn, {\em Singular elements of semi-simple algebraic groups}, Actes du Congr\`es International des Math\'ematiciens (Nice, 1970), 2, Gauthier-Villars, 279--284, (1971).

\item[\hbox to\parindent{\enskip\mbox{[DD]}}] T. De Fernex, R.
Docampo {\em Jacobian discrepancies and rational singularities},
arXiv math.AG/1106.2172v3.

\item[\hbox to\parindent{\enskip\mbox{[DEI]}}] T. De Fernex, L.
Ein and S. Ishii, {\em Divisorial valuations via arcs}, Publ. Res.
Inst. Math. Sci. 44 no. 2, 425--448 (2008).


\item[\hbox to\parindent{\enskip\mbox{[EM]}}] L. Ein and M.
Musta\cedilla{t}\v{a}, {\em Jet schemes and singularities}, Proc. Symp. Pure
Math. 80.2, 505--546 (2009).

\item[\hbox to\parindent{\enskip\mbox{[Hi1]}}] H. Hironaka, {\em
Resolution of singularities of an algebraic variety over a field
of characteristic zero: I and II}, Annals of Maths. 79 ns.1 and 2, 109--326
(1964).

\item[\hbox to\parindent{\enskip\mbox{[Hi2]}}] H. Hironaka, {\em
Characteristic polyhedra of singularities}, J. Math. Kyoto Univ. 7
n.3, 251--293 (1967).


\item[\hbox to\parindent{\enskip\mbox{[Is]}}] S. Ishii, {\em
Mather discrepancy and the arc spaces}, arXiv math.AG/1106.0345v1. 
to appear in Ann. de l'Institut Fourier, 62, (2012).

\item[\hbox to\parindent{\enskip\mbox{[Ka]}}] M. Kawakita, {\em Towards boundedness of minimal log discrepancies by Riemann-Roch Theorem},  Amer. J. Math. 133, no. 5, 1299--1311, (2011).

\item[\hbox to\parindent{\enskip\mbox{[KM]}}] J. Koll\'ar and S. Mori, {\em Birational 
Geometry of Algebraic Varieties}, Cambridge Tracts in Math., 134, (1998).

\item[\hbox to\parindent{\enskip\mbox{[KS]}}]  A. Kas  and M. Schlessinger, {\em
On the versal deformation of a complex space with an isolated singularity}, Math. Ann, 
196 23--29 (1972).

\item[\hbox to\parindent{\enskip\mbox{[Ko]}}] J. Koll\'ar,  {\em Rational Curves on
Algebraic Varieties}, Springer-Verlag, Ergebnisse der Math. 32, (1995).

\item[\hbox to\parindent{\enskip\mbox{[Mar]}}] D. Markushevich, {\em Minimal discrepancy for a cDV singularity is 1}, J. Math. Sci. Tokyo, 3 445--446, (1996).

\item[\hbox to\parindent{\enskip\mbox{[M]}}] H. Matsumura, {\em Commutative Ring 
Theory}, Cambridge studies in Advanced Math. 8, (1986).

\item[\hbox to\parindent{\enskip\mbox{[Re]}}] A.J. Reguera, {\em
Towards the singular locus of the space of arcs}, Amer. J. Math.
131, 2, 313-350, (2009).

\item[\hbox to\parindent{\enskip\mbox{[Sh]}}] V. V. Shokurov, {\em Letters of a bi-rationalist IV. Geometry of log flips},
   preprint (2002), arXiv:math/0206004.
\\


\end{itemize}

\noindent Shihoko Ishii, \\ Graduate School of Mathematical Science,
University of Tokyo, \\ 
3-8-1  Komaba, Meguro, 153-8914 Tokyo, Japan. \\
E-mail: shihoko@ms.u-tokyo.ac.jp \\

\noindent Ana J. Reguera, \\ Dep. de \'Algebra, Geometr\'ia y Topolog\'ia,  Universidad de Valladolid,\\
Paseo Bel\'en 7,
47011 Valladolid, Spain. \\
E-mail: areguera@agt.uva.es

\end{document}